\pdfoutput=1
\documentclass{article}

\usepackage[utf8]{inputenc}
\usepackage[english]{babel}
\usepackage{amsfonts}
\usepackage{amssymb}
\usepackage{mathrsfs}
\usepackage{amsthm}
\usepackage{latexsym,amsmath,amstext,amssymb,mathtools}
\usepackage{caption}
\usepackage{commath} 
\usepackage{hhline}
\usepackage{listings}
\usepackage{enumerate}
\usepackage{color}
\usepackage{xfrac}
\usepackage{faktor}
\usepackage{colonequals}
\usepackage{tikz-cd}
\usepackage{xr} 
\usepackage[toc,page]{appendix} 
\usepackage{hyphenat}
\usepackage{stmaryrd} 
\usepackage{csquotes}
\usepackage{array} 
\usepackage{indentfirst}
\usepackage{float}
\usepackage{capt-of}   

\usepackage[hidelinks]{hyperref}
\usepackage{cleveref}

\usepackage{tikz}
\usepackage{standalone}
\usepackage{forest}

\numberwithin{equation}{section}

\theoremstyle{plain}
\newtheorem{theorem}[equation]{Theorem}
\newtheorem*{theorem*}{Theorem}

\newtheorem{proposition}[equation]{Proposition}
\newtheorem{lemma}[equation]{Lemma}
\newtheorem*{lemma*}{Lemma}
\newtheorem*{proposition*}{Proposition}
\newtheorem*{corollary*}{Corollary}

\newtheorem{corollary}[equation]{Corollary}
\newtheorem{example}[equation]{Example}

\newtheorem{myremark}[equation]{Remark}
\newtheorem{maintheorem}{Theorem}

\newtheorem{introlemma}{Lemma}

\newtheorem{introproposition}{Proposition}
\newtheorem{conjecture}{Conjecture}
\newtheorem{mainexample}{Example}

\theoremstyle{definition}
\newtheorem{definition}[equation]{Definition}

\theoremstyle{remark}
\newtheorem*{remark}{Remark}

\newcommand{\dd}{\mathrm{d}}
\newcommand{\BO}{\mathcal{O}}

\newcommand\HH{\mathrm{H}}

\newcommand{\CH}{\mathrm{CH}}
\newcommand{\hh}{\mathrm{h}}

\newcommand{\CC}{\mathbb{C}}
\newcommand{\ZZ}{\mathbb{Z}}
\newcommand{\NN}{\mathbb{N}}
\newcommand{\QQ}{\mathbb{Q}}
\newcommand{\PP}{\mathbb{P}}

\newcommand\restr[1]{\raisebox{-.5ex}{$|$}_{#1}}
\newcommand{\sing}{\mathrm{sing}}
\newcommand{\Sing}{\mathrm{Sing}}

\newcommand{\sm}{\mathrm{sm}}
\newcommand{\codim}{\mathrm{codim}}

\newcommand{\GG}{\mathcal{G}}

\newcommand{\EE}{\mathscr{E}}
\newcommand{\FF}{\mathscr{F}}

\DeclareMathOperator{\divv}{div}

\newcommand{\sheafhom}{\mathscr{H}\kern -.5pt om}

\DeclareMathOperator{\Id}{Id}

\newcommand{\IC}{\mathrm{IC}}

\newcommand{\SL}{\mathrm{SL}}

\newcommand{\Nm}{\mathrm{Nm}}
\newcommand{\Pic}{\mathrm{Pic}}

\DeclareMathOperator{\Prym}{Prym}

\newcommand{\dec}{\mathrm{dec}}

\newcommand{\F}{\mathcal{F}}

\newcommand{\ud}{{\underline{d}}}

\newcommand{\SE}{\mathscr{S}}
\newcommand{\BE}{\mathscr{BE}}

\newcommand{\Ker}{\mathrm{Ker}}

\usepackage{microtype} 
\usepackage[backend=biber, style=alphabetic, giveninits=true, maxalphanames=99, doi=false,isbn=false,url=false]{biblatex}
\newcommand{\ns}{\mathfrak{X}}

\newcommand{\cc}{\mathrm{cc}}
\newcommand{\Xiasing}{\Xi_{\sing,\mathrm{ad}}}
\newcommand{\CL}{\mathfrak{L}}
\DeclareMathOperator{\Sl}{Sl}
\DeclareMathOperator{\Sp}{Sp}
\DeclareMathOperator{\Alt}{Alt}
\newcommand{\fJac}{\mathcal{J}^\mathrm{fake}}

\DeclareMathOperator{\Rep}{Rep}
\newcommand{\bb}{b}
\DeclareMathOperator{\Supp}{Supp}

\title{Singularities of theta divisors and the Tannakian Schottky problem}
\author{Constantin Podelski}

\begin{document}

\maketitle

\begin{abstract}
  \,With the Tannakian formalism, one can attach to any principally polarized abelian variety a reductive group, along with a representation. We obtain a completely new characterization of Jacobians relying on this Tannakian data and certain characteristic classes for singular varieties. As a corollary, we show that the Tannakian Schottky conjecture holds in dimension up to $5$. More generally, we show that the Tannakian Schottky conjecture holds on the bielliptic Prym locus in all dimensions. This provides the first non-trivial cases of this conjecture.
\end{abstract}



\section{Introduction}
\label{sec:intro}

Let $\mathcal{A}_g$ denote the moduli space of principally polarized abelian varieties over the complex numbers (ppav's for short) of dimension $g$. Denote by $\mathcal{J}_g\subset \mathcal{A}_g$ the locus of Jacobians. It is well known that for $g\geq 4$, the above inclusion is strict, and the Schottky problem is to characterize the Jacobian locus $\mathcal{J}_g$ inside $\mathcal{A}_g$. This problem has a very long history, with many successful approaches, as well as a few open conjectures, and we refer to Grushevsky's excellent paper for a survey \cite{Gru10SchottkyProb}. The singularities of theta divisors play a central role in this story. It is well known that the theta divisor of a non-hyperelliptic Jacobian of dimension $g$ has singular locus of dimension $g-4$. A famous result of Andreotti and Mayer \cite{Andreotti1967} states that the Jacobian locus $\mathcal{J}_g$ is an irreducible component of the Andreotti-Mayer locus
\[ \mathcal{N}^g \coloneqq \{ (A,\Theta)\in \mathcal{A}_g \,|\, \dim \Sing(\Theta)\geq g-4 \} \subset \mathcal{A}_g \,. \]
But this is only a weak solution to the Schottky problem. For $g\geq 4$, the Andreotti-Mayer locus has other irreducible components besides $\mathcal{J}_g$. One possible refinement is to look at characteristic classes of singular varieties. For a (possibly singular) subvariety $Z\subset A$ of an abelian variety $A$, we define a characteristic class denoted $c_\ast(Z)\in \HH_\ast(A)\coloneqq \HH_\ast(A,\ZZ)$ (\ref{Equ: Intro: Def of characteristic class of subvariety}). If $Z$ is smooth, and not invariant under translations by some positive-dimensional abelian subvariety, we recover the usual (signed) characteristic class
\[ c_\ast(Z)=c_{\ast}(T^\vee Z)\cap [Z] \,. \]
For an irreducible ample divisor $\Theta \subset A$, the dimension of the singular locus can be read off $c_\ast(\Theta)$ (Proposition \ref{Prop: Basic properties Chern Mather classes})
\[ \dim \Sing(\Theta)=\max \{k\,|\, c_k(\Theta)\neq [\Theta]^{g-k} \} \,. \]
We compute this class in the case of Jacobians (Proposition \ref{Prop: Chern-Mather classes of IC for Jacobians}):
\begin{introproposition}\label{introprop: chern class Jacobians}
    If $(JC,\Theta)\in\mathcal{J}_g$ is the Jacobian of a smooth curve $C$, then
\begin{equation}\label{Equ: Intro: Chern-Mather Jacobians} c_\ast(\Theta)=\begin{cases}
       \sum_{k=1}^g \bb_{k-1}[\Theta]^{ k}/k!\,, & \text{if $C$ non-hyperelliptic,} \\
   \sum_{k=1}^g \left(\bb_{k-1} -\bb^{2k-2}_{k-2}\right)[\Theta]^{k}/k!\,, & \text{if $C$ hyperelliptic,}
\end{cases}
\end{equation}
where $b_k\coloneqq \binom{2k}{k}$ and $b^n_k\coloneqq \binom{n}{k}$.
\end{introproposition}

It turns out that Proposition \ref{introprop: chern class Jacobians} is still not enough to characterize Jacobians, because for $g\geq 5$, the locus $\mathcal{N}^g$ has an irreducible component, called $\EE_{g,0}$ in \cite{Debarre1988}, for which the following is true (Lemma \ref{Lem: Chern Mather classes for E'g0}).
\begin{introproposition}
    If $g\geq 4$ and $(A,\Theta)\in \EE_{g,0}$ is general, then
    \[ c_\ast(\Theta)=\sum_{k=1}^g \bb_{k-1}[\Theta]^{ k}/k! \,. \]
\end{introproposition}
\textbf{The Tannakian framework.} To characterize Jacobians we need one more piece of (representation-theoretic) data, which is given by the Tannakian formalism. Let us briefly recall the setup (we refer to \cite{Weissauer2008TannCatAbVar, KraemerWeissauer2015VanThm} for a full account). Let $P(A)$ denote the category of perverse sheaves on $A$, and $\overline{P}(A)$ the quotient of $P(A)$ by the subcategory of negligible sheaves, i.e. those with vanishing Euler characteristic. The addition map $a:A\times A\to A$ induces a convolution product on $\overline{P}(A)$ denoted by $\ast$. For a subvariety $Z\subset A$, denote by $\IC_Z$ the intersection complex of $Z$. Let $\langle \IC_Z \rangle \subset \overline{P}(A)$ denote the full subcategory of $\overline{P}(A)$ containing all subquotients of convolution powers of $\IC_Z\oplus \IC_Z^\vee $. By \cite[Th. 13.2]{KraemerWeissauer2015VanThm}, this is a Tannakian category, i.e. there is a reductive group $G_Z$ over $\CC$ and an equivalence of categories
\[ \omega: \langle \IC_Z \rangle \overset{\sim}{\longrightarrow} \Rep(G_Z)\,,\]
where $\Rep(G_Z)$ is the category of finite-dimensional representations of $G_Z$. We denote by $\omega_Z\coloneqq \omega(\IC_Z)$ the representation corresponding to $Z$ under this equivalence. We say that $G_Z$ (resp. $\omega_Z$) is the Tannakian group (resp. representation) associated to $Z\subset A$. Under the above equivalence we have \cite[Rem. 1.1]{Kraemer2015}
\begin{equation}\label{IntroEqu: dim omega(P) is euler characteristic of P} 
\dim \omega_Z=\chi(\IC_Z)\coloneqq \sum_k (-1)^k \hh^k(A,\IC_Z)\,. 
\end{equation}
Given a ppav $(A,\Theta)$, the group $G_\Theta$ does depend on the chosen translate of $\Theta$, but $\tilde{G}_\Theta \coloneqq [ G_\Theta^\circ,G_\Theta^\circ]$ (the derived group of the connected component containing the identity) does not \cite[Lem. 4.3.2]{Kraemer2021MicrolocalGauss2}. From now on, we will always implicitly fix a symmetric translate of $\Theta$. \par 
Intuitively, one can think of $(G_\Theta,\omega_\Theta)$ as a categorification of the monodromy of the Gauss map of $\Theta$ (this statement is made precise in \cite[Th. 1.7]{Kraemer2022MicrolocalGauss1}). In the case of Jacobians, this group and representation take a very specific form \cite{Weissauer2006BrillNoetherSheaves}:
\begin{mainexample}\label{mainexample: Tannakian group of Jacobian}
    Let $(JC,\Theta)$ be the Jacobian of a smooth curve of genus $g$. For $n=g-1$, we have
    \begin{equation}\label{Equ: Tannakian Group and Rep for Jacobians} 
        G_\Theta=\begin{cases}
            \Sl_{2n}(\CC) / \mu_{n} \\
            \Sp_{2n}(\CC)/ \mu_{\gcd(2,n)}
        \end{cases} \quad \text{and}\qquad \omega_\Theta=\begin{cases}
            \Alt^{n}(\CC^{2n}) \\
            \Alt^{n}(\CC^{2n})/\Alt^{n-2}(\CC^{2n})
        \end{cases}
    \end{equation}
    in the non-hyperelliptic, resp. hyperelliptic case (for some symmetric translate of $\Theta$).
\end{mainexample}

To state our main result, we need a technical condition involving the \emph{problematic locus} $S_\Theta \subset T^\vee_0 A$. We defer its definition until later in the introduction (\ref{Equ: Intro: Problematic locus}). Essentially, it is a finiteness assumption on the Gauss map of $\Theta$. Our main result is that up to this finiteness assumption on the Gauss map of $\Theta$, the characteristic class $c_\ast(\Theta)$ together with the Tannakian formalism characterizes Jacobians.
\begin{maintheorem}\label{maintheorem: Tannakian criterion for Jacobians}
    Let $(A,\Theta)\in \mathcal{A}_g$ be such that 
\begin{itemize}
        \item $(G_\Theta,\omega_\Theta)$ is as in (\ref{Equ: Tannakian Group and Rep for Jacobians}) for a suitable choice of symmetric translate $\Theta$.
        \item For $r\in \{1,2\}$ we have
        \begin{align*}
    c_r(\Theta)= \begin{cases}
        \bb_{g-r-1}[\Theta]^{ g-r}/(g-r)! & \text{(non-hyperelliptic case)} \\
    \left(\bb_{g-r-1} -\bb^{2g-2r-2}_{g-r-3}\right)[\Theta]^{g-r}/(g-r)! & \text{(hyperelliptic case).}
    \end{cases} 
    \end{align*}
    \item We have $\codim_{T^\vee_0 A} S_\Theta > 2\,.$
\end{itemize}
Then $(A,\Theta)$ is a Jacobian.
\end{maintheorem}
We prove Theorem \ref{maintheorem: Tannakian criterion for Jacobians} in Section \ref{Subsec: proof of maintheorem Tannakian criterion}. The main idea of the proof is to use relations on Chern classes coming from representation theory and $\lambda$-rings to show that if the conditions of the above theorem hold, then $\Theta$ is a sum of curves. This implies that $(A,\Theta)$ is a Jacobian by a theorem of Schreieder \cite{Schreieder2015ThetaCurveSummand}.

Following \cite{Kraemer2021MicrolocalGauss2}, we say that a ppav $(A,\Theta)\in\mathcal{A}_g$ is a \emph{fake Jacobian} (resp. non-hyperelliptic fake Jacobian, hyperelliptic fake Jacobian) if for some translate of $\Theta$, the group and representation $(G_\Theta,\omega_\Theta)$ take the form of Example \ref{mainexample: Tannakian group of Jacobian} (resp. the first or the second form of the example). In the spirit of the Andreotti-Mayer locus, Kr{\"a}mer defines the locus of fake Jacobians by 
\[ \fJac_g \coloneqq  \overline{\left\{ (A,\Theta)\in \mathcal{A}_g \,|\, (G_\Theta,\omega_\Theta) \text{ as in (\ref{Equ: Tannakian Group and Rep for Jacobians})}  \right\}} \subset \mathcal{A}_g\,. \] 
By Kr{\"a}mer \cite{Kraemer2021MicrolocalGauss2}, this gives a weak solution to the Schottky problem, i.e. the Jacobian locus $\mathcal{J}_g$ is an irreducible component of $\fJac_g$.
Kr{\"a}mer and Weissauer make the following conjecture:
\begin{conjecture}[Tannakian Schottky]\label{Conj: Tannakian Schottky}
    For $g\geq 0$, we have
    \[ \fJac_g=\mathcal{J}_g\,. \]
\end{conjecture}
 For $g\leq 3 $ we have $\mathcal{J}_g=\mathcal{A}_g$, thus Conjecture \ref{Conj: Tannakian Schottky} is immediate. As a consequence of Theorem \ref{maintheorem: Tannakian criterion for Jacobians} we prove the first two non-trivial cases of Conjecture \ref{Conj: Tannakian Schottky}, namely $g=4$ and $g=5$.
\begin{maintheorem}\label{maintheorem: Tannakian Schottky in dim up to 5}
    For $0\leq g \leq 5$, we have
    \[ \fJac_g=\mathcal{J}_g\,. \]
\end{maintheorem}
We will prove this in Section \ref{Sec: The dimension 5 case}. Let us examine the case $g=5$ more closely. By \cite[Lem. 4.2.1]{Kraemer2021MicrolocalGauss2}, we have
\[  \fJac_5\subset  \mathcal{N}^5 \,. \] 
Recall that a curve is \emph{bielliptic} if it admits a double cover to an elliptic curve. Denote by $\BE_g\subset \mathcal{A}_g$ the locus of \emph{bielliptic Prym varieties}, i.e. Pryms arising from étale double covers of bielliptic curves (see Section \ref{Sec: The bielliptic Prym locus}). By Donagi \cite{Donagi1981tetragonal} and Debarre \cite{Debarre1988}, the irreducible components of $\mathcal{N}^5$ are given by
\begin{equation}\label{equ: Andreotti dim 5}
    \mathcal{N}^{5} = \mathcal{J}_5 \cup \mathcal{A}_{1,4} \cup \EE_{5,0} \cup \EE_{5,1} \cup \EE_{5,2} \,,
\end{equation} 
where $\mathcal{A}_{1,4}$ is the locus of ppav's that are product of an elliptic curve with a $4$-dimensional ppav, and $\EE_{5,0},\EE_{5,1},\EE_{5,2}$ are the three irreducible components of the bielliptic Prym locus $\BE_5$. One can immediately exclude fake Jacobians in $\mathcal{A}_{1,4}$, since for decomposable ppav's the group $G_\Theta$ is trivial. Theorem \ref{maintheorem: Tannakian Schottky in dim up to 5} then follows from the general fact that the Tannakian formalism characterizes (non-hyperelliptic) Jacobians on the bielliptic Prym locus, in all dimensions:
\begin{maintheorem}\label{maintheorem: Tannakian Schottky on Bielliptic Prym locus}
Let $g\geq 0 $ and $(A,\Theta)\in \BE_g$ be a non-hyperelliptic fake Jacobian. Then $(A,\Theta)$ is a (non-hyperelliptic) Jacobian.
\end{maintheorem}
We prove Theorem \ref{maintheorem: Tannakian Schottky on Bielliptic Prym locus} in Section \ref{Sec: The proof of Theorem 3}. The idea of the proof is the following. For $(A,\Theta)\in\mathcal{A}_g$, we have $\dim \omega_\Theta = \deg c_\ast(\Theta)$ by (\ref{IntroEqu: dim omega(P) is euler characteristic of P}). Surprisingly, it turns out that for all $(P,\Xi)\in \BE_g$ such that $\deg c_\ast(\Theta)=b_{g-1}$, the ppav $(P,\Xi)$ also satisfies conditions ii) and iii) of Theorem \ref{maintheorem: Tannakian criterion for Jacobians}. Theorem \ref{maintheorem: Tannakian Schottky on Bielliptic Prym locus} thus follows from Theorem \ref{maintheorem: Tannakian criterion for Jacobians}. \par

\textbf{The problematic locus.} For a subvariety $Z\subset A$, the \emph{conormal variety} is the Zariski closure of the conormal bundle to the smooth locus $Z_\sm$, taken inside the total space of the cotangent bundle:
\[ \Lambda_Z \coloneqq \overline{ N^\vee_{Z_\sm}A} \subset T^\vee A \simeq A\times V \,, \]
where $V\coloneqq T^\vee_0 A$. 
The group of \emph{Lagrangian cycles} $\mathscr{L}(A)$ is the free abelian group generated by conormal varieties to subvarieties of $A$. The \emph{Chern-Mather class} of $\Lambda_Z$ is a class $c_M(\Lambda_Z)\in \HH_\ast(A)$ (see (\ref{Equ: Definition of Chern-Mather class}) for a definition). 
We say that $\Lambda_Z$ is \emph{negligible} if $\deg c_M(\Lambda_Z)=0$, and we call it \emph{clean} otherwise. By \cite{Weissauer2011DegenerateGauss}, $\Lambda_Z$ is negligible if and only if $Z$ is invariant under translation by a positive-dimensional abelian subvariety of $A$. Denote by $\mathfrak{L}(A)\subset \mathscr{L}(A)$ the subgroup generated by clean conormal varieties. If $K\in P(A)$ is a perverse sheaf, then its characteristic cycle is
\[ \mathrm{CC}(K)=\sum_i n_i \Lambda_{Z_i} \in \mathscr{L}(A) \]
for some subvarieties $Z_i\subset A$ and $n_i\in \NN_{>0}$. We then define the \emph{clean characteristic cycle} as the clean part in the above sum
\[ \cc(K)\coloneqq \sum_{i\,|\, \Lambda_{Z_i} \text{ clean}} n_i \Lambda_{Z_i} \,.\]
By Kashiwara's index formula, a perverse sheaf $K\in P(A)$ is negligible if and only if its characteristic cycle is negligible. We define the (clean) Chern class on $\overline{P}(A)$ by
\[ c_\ast \coloneqq c_M\circ \mathrm{cc}: \overline{P}(A) \to \HH_\ast(A) \,. \] 
For a subvariety $Z\subset A$, we define its (clean) characteristic class by
\begin{equation}\label{Equ: Intro: Def of characteristic class of subvariety}
    c_\ast(Z)\coloneqq c_\ast(\IC_Z)\,.
\end{equation} 
The addition map on $A$ induces a natural ring structure on $\mathfrak{L}(A)$ whose product we denote by $\circ$ (see Section \ref{Subsec: The ring of Lagrangian cycles}). In this setting, the clean characteristic cycle
\[ \mathrm{cc}:\left(\overline{P}(A),\ast \right) \to \left( \mathfrak{L}(A),\circ \right) \] 
is a ring morphism. Consider the ring structure on $\HH_\ast(A)$ induced by the Pontryagin product. The Chern-Mather class induces a map $c_{\ast}: \mathfrak{L}(A)\to \HH_\ast(A)$, but this is in general not a ring morphism. We want to understand when $c_\ast : \langle \IC_Z \rangle \to \HH_\ast(A)$ is a ring morphism. Let $\Lambda\coloneqq \mathrm{cc}(\IC_Z)$. We stratify $ V$ by the locally closed sets
\[ S_{\Lambda,d}'\coloneqq \{ \xi \in V\setminus \{0\} \,|\, \dim \langle \Lambda\restr{\xi} \rangle = d \}\subset  V \,, \quad \text{for $d\geq 0$,} \] 
where for a subset $W\subset A$, we define $\langle W \rangle$ to be the smallest abelian subvariety of $A$ containing the connected component of $0$ of $W-W$. We define the \emph{problematic locus} as
\begin{equation}\label{Equ: Intro: Problematic locus}
S_Z \coloneqq \bigcup_{d=1}^g  \bigcup_{ \substack{S\subseteq  S'_{\Lambda,d} \\ \codim_{V} S \leq d }} \overline{S}         \subset V \,. 
\end{equation}
Consider the quotient 
\[ \HH_{\leq d}(A)\coloneqq \faktor{\HH_{\ast}(A)}{\HH_{>d}(A)} \,, \quad \text{with} \quad \HH_{>d}(A) \coloneqq \bigoplus_{k> d} \HH_{k}(A)\,. \] 
We have the following (Lemma \ref{Lem: Chern-Mather: multiplicativity relative to codim of problematic locus}):
\begin{introlemma}
     For $Z\subset A$ and $d\geq 0$, the Chern class in degree up to $d$
    \[ c_{\leq d}:\left(\langle \IC_Z \rangle,\ast \right) \to \left(\HH_{\leq d}(A) ,\ast \right)\]
    is a ring morphism if and only if $d<\codim_V S_Z$.
\end{introlemma}
Being a Tannakian category, $\langle \IC_Z \rangle$ carries a natural $\lambda$-ring structure. $\HH_\ast(A)$ carries a $\lambda$-ring structure as well, whose Adams operations $\psi^n$ correspond to $[n]_\ast:\HH_\ast(A)\to \HH_\ast(A)$. As a consequence of \cite[Prop. 1.3.3]{Kraemer2021MicrolocalGauss2}, whenever $c_\ast$ is a map of rings, it is automatically a map of $\lambda$-rings. Thus, there is a lot of structure to work with. This is one of the main ingredients in the proof of Theorem \ref{maintheorem: Tannakian criterion for Jacobians}.
\begin{remark}\,
\begin{enumerate}
    \item We have $\codim_V S_Z \geq 2$ (Remark \ref{Rem: S Lambda 1 is empty}).
    \item The problematic locus is linear, i.e. 
    \[ S_Z=\bigcup_{i=1}^k T^\vee_0 B_i \subset T^\vee_0 A \,, \]
    where $B_1,\dots,B_k$ are abelian subvarieties of $A$ (Lemma \ref{Cor: problematic locus is linear}). In particular, $S_Z$ is empty if $A$ is simple.
    \item The problematic locus $S_Z$ is contained in the locus where the Gauss map $\gamma_Z:\mathrm{cc}(\IC_Z)\to V$ has positive-dimensional fibers. In particular, it is empty if $\gamma_Z$ is finite.
    \item If $Z\subset A$ is a smooth subvariety with ample normal bundle, the Gauss map $\gamma_{Z}$ is finite \cite[Prop. 1.1]{Debarre1995FultonHansen}, so $S_{Z}$ is empty.
    \item The condition of having an empty problematic locus is strictly stronger than being geometrically non-degenerate in the sense of \cite{ran1980subvarieties,Debarre1995FultonHansen} (Proposition \ref{Prop: empty problematic locus implies geometrically non-degenerate} and Example \ref{Ex: non-empty problematic locus}).
\end{enumerate} 
\end{remark}

As illustrated by the above remarks, the failure of the Chern class to be multiplicative on $\langle \IC_Z \rangle$ is intimately related to the singularities of the subvariety $Z\subset A$. This seems to be one of the challenges when dealing with singular varieties. For example, most recent applications of the Tannakian formalism \cite{javanpeykar2022monodromy, kraemer2024arithmeticfinitenessirregularvarieties, kraemer2024tannakagroupe6arises} only deal with smooth subvarieties with ample normal bundle in $A$. One might hope that subvarieties with empty problematic locus behave similarly (from a Tannakian perspective) as smooth subvarieties with ample normal bundle.
\par

\par 
This text is organized as follows. In Section \ref{Sec: Tannakian Schottky} we recall the Tannakian formalism on abelian varieties, study the ring of Lagrangian cycles and the multiplicativity of the Chern-Mather class. In Section \ref{A criterion for detecting Jacobians} we prove Theorem \ref{maintheorem: Tannakian criterion for Jacobians}. In Section \ref{Sec: The bielliptic Prym locus} we recall the main properties of the bielliptic Prym locus $\BE_g$. We compute the characteristic cycle of the corresponding Prym theta divisors (Theorem \ref{Theorem: Characteristic Cycle for E'gt}), study the fibers of the Gauss map (Lemma \ref{Lem: Gauss map finite for Eg0'}) and compute the Chern-Mather class (Lemma \ref{Lem: Chern Mather classes for E'g0}). As a result, we obtain that bielliptic Pryms satisfy the conditions of Theorem \ref{maintheorem: Tannakian criterion for Jacobians}. In Section \ref{Sec: The proof of Theorem 3}, we put all of this together and prove Theorem \ref{maintheorem: Tannakian Schottky on Bielliptic Prym locus}. Finally, we prove Theorem \ref{maintheorem: Tannakian Schottky in dim up to 5} in Section \ref{Sec: The dimension 5 case} by excluding the hyperelliptic case in dimension up to $5$ by a divisibility argument on the Chern-Mather class.

\par 
\textbf{Acknowledgments.}{The author would like to thank his supervisor Prof. Thomas Kr{\"a}mer for his insights and the many discussions that led to this paper.}

\section{The Tannakian formalism and some conormal geometry}\label{Sec: Tannakian Schottky}

\subsection{The Tannakian formalism}\label{Sec: The Tannakian formalism}
We briefly recall the construction of the Tannakian category. We refer to \cite{KraemerWeissauer2015VanThm} for a detailed account and proofs. Let $(A,\Theta)\in\mathcal{A}_g$ be a ppav and $D(A)$ denote the derived category of bounded complexes of constructible sheaves with $\CC$-coefficients on $A$. Let $S(A)\subset D(A)$ denote the full category of \emph{negligible} complexes, i.e. complexes $K\in D(A)$ satisfying
\[\chi(\,^p \mathcal{H}^n(K))=0 \,, \quad \text{for all $n\in \ZZ$,}\]
where $\,^p\mathcal{H}^n(K)$ is the $n$-th perverse cohomology of $K$ and $\chi$ is the Euler characteristic is defined by
\[ \chi(K)\coloneqq \sum_n (-1)^n \dim \HH^n(A,K)\,. \]
Negligible perverse sheaves correspond to perverse sheaves whose simple subquotients are invariant under translation by some abelian subvariety of positive dimension. $S(A)$ is a thick subcategory of $D(A)$ and one defines $\overline{D}(A)$ as the triangulated quotient category. Let $P(A)\subset D(A)$ denote the category of perverse sheaves and $\overline{P}(A)\subset \overline{D}(A)$ its image in the quotient category. The addition map $a:A\times A\to A$ induces a convolution product on $D(A)$
\begin{align*}
    \ast:D(A)\times D(A) &\to D(A) \,, \qquad K,L \mapsto Ra_\ast (K\boxtimes L) \,.
\end{align*}
This convolution product descends to $\overline{D}(A)$, preserves $\overline{P}(A)$ and induces on both a structure of rigid symmetric monoidal category. Let $Z\subset A$ be a subvariety of $A$. We denote by $\IC_Z\in P(A)$ the \emph{intersection cohomology complex} associated to $Z$, i.e. the unique simple perverse sheaf on $Z$ that restricts to $\underline{\CC}_{Z_\sm}[d]$ on the smooth locus $Z_\sm$, where $d=\dim Z$ (see \cite[Ex. 8.4.7]{maxim2019}). Let $\langle \IC_Z \rangle \subset \overline{P}(A)$ denote the full subcategory of $\overline{P}(A)$ containing all quotients of convolution powers of $\IC_Z\oplus \IC_Z^\vee $. This is a rigid symmetric monoidal category with respect to the convolution \cite[Th. 13.2]{KraemerWeissauer2015VanThm}. Thus, there is a reductive group $G_Z$ over $\CC$ and an equivalence
\[ \omega: \langle \IC_Z \rangle \overset{\sim}{\longrightarrow} \Rep_\CC(G_Z)\]
with the rigid symmetric monoidal category of finite-dimensional algebraic representations of $G_Z$. We denote by $\omega_Z\coloneqq \omega(\IC_Z)$ the representation corresponding to $Z$ under this equivalence. We say that $G_Z$ (resp. $\omega_Z$) is the Tannakian group (resp. representation) associated to $Z\subset A$. Under the above equivalence we have
\begin{equation}\label{Equ: dim omega(P) is euler characteristic of P} 
\dim \omega_Z=\chi(\IC_\Theta)\,. 
\end{equation}
Given a ppav $(A,\Theta)$, the group $G_\Theta$ does depend on the chosen translate of $\Theta$, but the derived group of the connected component containing the identity $\tilde{G}_\Theta\coloneqq [G_\Theta^\circ,G_\Theta^\circ]$ does not \cite[Lem. 4.3.2]{Kraemer2021MicrolocalGauss2}. From now on, we will always implicitly fix a symmetric translate of $\Theta$. 
\par

\subsection{The ring of Lagrangian cycles}\label{Subsec: The ring of Lagrangian cycles}
Let $A$ be a $g$-dimensional abelian variety, and $V\coloneqq T^\vee_0 A$. The cotangent bundle $T^\vee A = A\times V$ is canonically trivialized using translations. Let $\Lambda \subset A\times V$ be a $g$-dimensional subvariety which is \emph{conic}, i.e. invariant under rescaling of the fibers of the cotangent bundle. We define $\PP\Lambda\subset A\times \PP V$ to be the projectivization. We define the \emph{Gauss map} on $\Lambda$ and $\PP\Lambda$ to be the projection onto the second factor
\[ \gamma_\Lambda : \Lambda \to V \,, \quad \gamma_{\PP\Lambda}: \PP\Lambda \to \PP V \,. \]
The \emph{degree} of $\Lambda$ is $\deg(\Lambda)\coloneqq \deg \gamma_\Lambda$. We say that $\Lambda$ is \emph{negligible} if $\deg(\Lambda)=0$, i.e. if $\gamma_\Lambda$ is not dominant, and \emph{clean} otherwise. For a subvariety $Z\subset A$ we define its \emph{conormal variety} by
\[ \Lambda_Z \coloneqq \overline{ \{ (x,\xi)\in T^\vee A \,|\, x\in Z_\sm \,, T_x Z \subset \Ker(\xi) \}} \subset A\times V \,, \]
and denote by $\gamma_Z\coloneqq \gamma_{\PP\Lambda_Z}:\PP\Lambda_Z\to \PP V$ the associated projectivized Gauss map. We say that $Z$ is negligible (resp. clean) if $\Lambda_Z$ is. Note that $Z$ is negligible if and only if $Z$ is invariant under translation by a positive-dimensional abelian subvariety $B\subset A$ \cite{Weissauer2011DegenerateGauss}. Let $B\subset A$ be the maximal abelian subvariety of $A$ such that $Z$ is invariant under translations by $B$, and let $q:A\to A/B$ be the projection. Then
\begin{equation}\label{Equ: negligible subvarieties}
Z=q^{-1}\left(W\right)\,, 
\end{equation}
where $W\coloneqq q(Z)$ is a clean subvariety of $A/B$. \par 
The group of \emph{Lagrangian cycles} $\mathscr{L}(A)$ is defined as the free abelian group generated by the conormal varieties $\Lambda_Z$ for subvarieties $Z\subset A$. The group of \emph{clean Lagrangian cycles} is defined as the subgroup of $\mathfrak{L}(A)\subset \mathscr{L}(A)$ generated by those conormal varieties
$\Lambda_Z$ that are clean. 
Following \cite{Kraemer2021MicrolocalGauss2}, we can endow $\mathfrak{L}(A)$ with a structure of $\lambda$-ring. Denote by $a:A\times A\to A$ the addition map and by $\delta:V\hookrightarrow V\times V$ the diagonal embedding. We define $\rho \coloneqq \Id_A\times \delta$, $\rho'\coloneqq \Id_{A^2}\times \delta$, $\varpi \coloneqq a\times \Id_{V^2}$ and $\varpi'\coloneqq a\times \Id_V$. We obtain the fiber square
\begin{equation}\label{Equ: Def or rho and varpi} 
\begin{tikzcd}
     A^2\times V \arrow[r,hook, "\rho'"]  \arrow[d,"\varpi'"'] & A^2 \times V^2  \arrow[d,"\varpi"] \\
     A\times V \arrow[r,hook,"\rho"']  & A\times V^2 \,.
\end{tikzcd}
\end{equation}
For $\Lambda_1,\Lambda_2\in \mathscr{L}(A)$, we define the \emph{non-clean product} by
\begin{equation}\label{Equ: Definition of alternate product of lagrangians}
\Lambda_1 \tilde{\circ} \Lambda_2 \coloneqq \rho^{!}\left(\varpi_\ast(\Lambda_1\times \Lambda_2)\right)\in  \CH_g(T^\vee A)\,. 
\end{equation}
Note that by the projection formula \cite[Th. 6.2]{Fulton1998} we have $\Lambda_1 \tilde{\circ} \Lambda_2= \varpi'_\ast\left(\rho'^{!}(\Lambda_1\times \Lambda_2)\right)$. By the proposition below, the locus $\rho^{-1}\left(\varpi (\Supp(\Lambda_1\times \Lambda_2))\right)$ is of dimension at most $g$ (and always of dimension $g$ when $\Lambda_1,\Lambda_2$ are clean). Thus, from now on we will view $\Lambda_1 \tilde{\circ} \Lambda_2$ as a cycle supported on $\rho^{-1}\left(\varpi (\Supp(\Lambda_1\times \Lambda_2))\right)$.
\begin{proposition}
    For $\Lambda_1,\Lambda_2\in\mathscr{L}(A)$, the non-clean product $\Lambda_1 \tilde{\circ} \Lambda_2\in \mathscr{L}(A)$ is a (possibly not clean) Lagrangian cycle.
\end{proposition}
\begin{proof}
Without loss of generality, we can assume that $\Lambda_1,\Lambda_2\subset T^\vee A$ are two conic Lagrangian subvarieties in $T^\vee A$. If $\dim \varpi(\Lambda_1\times \Lambda_2)<2g$, then $\Lambda_1\tilde{\circ} \Lambda_2=0$ (note this can happen only if either $\Lambda_1$ or $\Lambda_2$ is negligible). From now on we assume that $\dim \varpi(\Lambda_1\times \Lambda_2)=2g$. As $\rho$ is a regular embedding of codimension $g$, every irreducible component of $\rho^{-1}\left(\varpi (\Supp(\Lambda_1\times \Lambda_2))\right)$ is of dimension at least $g$. Denote by $\alpha\in \HH^0(T^\vee A, \Omega^1_{T^\vee A})$ the canonical $1$-form of $T^\vee A$ (see \cite[Sec. E.3]{hotta_d-modules_2008}). Recall that if $dx_1,\dots,dx_g$ denotes a basis of $T^\vee_0 A$, and $\xi_1,\dots,\xi_g$ denote the corresponding coordinates on $V=T^\vee_0 A$, then
\[ \alpha=\sum_{i=1}^g \xi_i d x_i \in \HH^0(T^\vee A, \Omega^1_{T^\vee A}) \,. \]
Let $\alpha_1\coloneqq p_1^\ast \alpha$ and $\alpha_2\coloneqq p_2^\ast \alpha$, where $p_1,p_2:(T^\vee A)^2\to T^\vee A$ denote the projection onto the first and second factor respectively. Denote by $d x_1'\dots,d x_g'$ the corresponding basis of the second factor of $(T^\vee_0 A)^2$ and by $\xi'_1,\dots,\xi'_g$ the corresponding coordinates on the second factor of $V^2$. We have
\[ \rho'^\ast(\alpha_1+\alpha_2)= \sum_{i=1}^g \xi_i (d x_i+ d x_i')= \varpi'^\ast \alpha \in \HH^0(A^2\times V,\Omega^1_{A^2\times V}) \,.\]
Recall that a conic subvariety $\Lambda \subset T^\vee A$ is isotropic if and only if $\alpha\restr{\Lambda}=0$ \cite[Lem. E.3.1]{hotta_d-modules_2008}. In particular, $(\alpha_1+\alpha_2)\restr{\Lambda_1\times\Lambda_2}=0$. Thus
\[\rho'^\ast(\alpha_1+\alpha_2)\restr{\rho'^{-1}(\Lambda_1\times \Lambda_2)}=0 \,,\]
thus $\alpha\restr{\Lambda_1 \tilde{\circ} \Lambda_2}=0$ and $\Lambda_1 \tilde{\circ} \Lambda_2$ is isotropic. Every irreducible component of $\Lambda_1 \tilde{\circ} \Lambda_2$ is of dimension at least $g$, thus $\Lambda_1 \tilde{\circ} \Lambda_2$ is Lagrangian.
\end{proof}
By the above proposition, if $\Lambda_1,\Lambda_2$ are effective cycles, we can write 
\begin{equation}\label{Equ: non-clean product is sum of clean and negligible with positive coefficients} \Lambda_1 \tilde{\circ} \Lambda_2 = \sum_{ Z\text{ clean}} n_Z \Lambda_Z+ \sum_{Z \text{ negligible}} n_Z \Lambda_Z \in \mathscr{L}(A)\,,
\end{equation}
where $n_Z\in \ZZ_{>0}$ by \cite[Prop. 7.1]{Fulton1998}. For $\Lambda_1,\Lambda_2\in\mathfrak{L}(A)$, we define the \emph{clean product} (or simply the product) $\Lambda_1\circ \Lambda_2\in \mathfrak{L}(A)$ as the clean part in the above sum. Clearly, this coincides with Kr{\"a}mer's definition \cite{Kraemer2021MicrolocalGauss2}
\[ \Lambda_1 \circ \Lambda_2 \coloneqq \overline{  \varpi'_\ast\left(\rho'^{-1}(\Lambda_1\times \Lambda_2)\restr{U} \right) } \in \mathfrak{L}(A)\,, \]
where $U\subset V$ is an open set over which both $\gamma_{\Lambda_1}$ and $\gamma_{\Lambda_2}$ are finite. The product $\circ$ induces a ring structure on $\mathfrak{L}(A)$. By \cite[Prop. 1.3.3]{Kraemer2021MicrolocalGauss2}, $\mathfrak{L}(A)$ comes with a natural $\lambda$-ring structure whose Adams operations $\Psi^n$ are given by $[n]_\ast : \mathfrak{L}(A)\to \mathfrak{L}(A)$ where $[n]:A\times V\to A\times V$ is the multiplication by $n$ on the $A$ factor. For $\Lambda\in \mathfrak{L}(A)$, we denote by $\langle \Lambda \rangle \subset \mathfrak{L}(A)$ the smallest $\lambda$-subring of $\mathfrak{L}(A)$ containing $\Lambda$ that is stable under taking irreducible components of its members.

\subsection{The problematic locus}\label{Subsec: The problematic locus}
Let $\Lambda\in \mathfrak{L}(A)$ and $\Lambda_1,\Lambda_2\in\langle \Lambda \rangle$. We want to understand when the non-clean product $\Lambda_1 \tilde{\circ} \Lambda_2$ and the (clean) product $\Lambda_1 \circ \Lambda_2$ differ. For $\Lambda\in \mathscr{L}(A)$, let
\begin{align*}
     S_{\Lambda,d}'&\coloneqq \{ \xi \in V\setminus \{0\} \,|\, \dim \langle \Lambda\restr{\xi} \rangle = d \}\subset  V \,, \quad \text{for $d\geq 0$,}
     \end{align*}
We define the \emph{problematic stratum} (resp. \emph{problematic locus}) by
\[ 
 S_{\Lambda,d} \coloneqq \bigcup_{ \substack{S\subseteq  S'_{\Lambda,d} \\\codim_{V} S \leq d }} \overline{S}\subset V\, \qquad \left(\text{resp.}\quad 
S_\Lambda \coloneqq \bigcup_{d=1}^g S_{\Lambda,d}\subset V  \right)\,.
\]

\begin{myremark}\label{Rem: S Lambda 1 is empty}
    Note that $\Lambda$ has $1$-dimensional fibers above $S_{\Lambda,1}$. If $\codim_{V} S'_{\Lambda,1} \leq 1$, then for dimension reasons an irreducible component of $\Supp(\Lambda)$ sits above $S_1$, contradicting the cleanness of $\Lambda$. Thus, $S_{\Lambda,1}$ is empty. Similarly, by the results below, $S_{\Lambda,g}$ is empty.
\end{myremark}

\begin{lemma}\label{Lem: negligible components appear}
    If $\Lambda_1,\Lambda_2\in \langle \Lambda \rangle$, then any negligible component of $\Lambda_1\tilde{\circ} \Lambda_2$ lies above the problematic locus $S_\Lambda$. Reciprocally, for $d\geq 1$ and an irreducible component $S\subset S_{\Lambda,d}$, there exist $\Lambda_1,\Lambda_2\in \langle \Lambda \rangle$ such that a negligible component of $\Lambda_1 \tilde{\circ}\Lambda_2$ lies above $S$.
\end{lemma}
\begin{proof}
    Let $\Lambda_1,\Lambda_2\in \langle \Lambda \rangle$ and let $\Lambda_Z$ be a negligible component of $\Lambda_1\tilde{\circ} \Lambda_2$. Then by the previous section, $Z$ is invariant by translation by a positive-dimensional abelian subvariety $B\subset A$, and the quotient $W\subset A/B\eqcolon B'$ is a clean subvariety of $B'$. Let $\Lambda_{W} \subset T^\vee B'$ be the conormal variety to $W$. Then
\begin{equation}\label{Equ: Lambda Z of negligible Z} \Lambda_Z= B\times_{B'} \Lambda_{W} \subset A\times T^\vee_0 B' \subset A\times T^\vee_0 A \,. 
\end{equation}
In particular, $\Lambda_Z$ sits above $T^\vee_0 B'$ and for a general $\xi\in T^\vee_0 B'$
\[ \langle \gamma_Z^{-1}(\xi) \rangle = B \,. \]
Thus, $ T^\vee_0 B' \subset S_{\Lambda,d}$ with $d\coloneqq \dim B$. \par 
Conversely, let $S$ be an irreducible component of $S_{\Lambda,d}$ for some $d$. Let $\xi \in S$ be generic and $n> 1$ be minimal such that 
\[ \dim \left(\Lambda^{\tilde{\circ}n}\right)\restr{\xi}=d \,.\]
By assumption $\codim_V S \leq d$, thus there is an irreducible component (necessarily negligible) of $\Lambda^{\tilde{\circ} n}$ lying above $S$.
\end{proof}
We have the following:
\begin{corollary}\label{Cor: problematic locus is linear}
    Let $S_{\Lambda,d}$ be a non-empty problematic stratum for some $d\geq 1$. Then $S_{\Lambda,d}$ is a finite union of linear spaces
    \[ S_{\Lambda,d}= \bigcup_{i=1}^k T^\vee_0 B_i \,, \]
    where $B_1,\dots,B_k$ are abelian subvarieties of $A$ of codimension $d$. In particular, $\dim S_{\Lambda,d}=g-d$.
\end{corollary}
\begin{proof}
    By the above lemma, any irreducible component of $S_{\Lambda,d}$ is realized by a negligible component $\Lambda_Z$ for some $Z$. The assertion follows from (\ref{Equ: Lambda Z of negligible Z}) and the fact that $A$ has only finitely many abelian subvarieties \cite[Th. 5.3.7]{Birkenhake2004}.
\end{proof}
Recall from \cite{ran1980subvarieties} that a subvariety $Z\subset A$ is said to be \emph{geometrically non-degenerate} if for all abelian quotients $q:A\to B$, either $q(Z)=B$ or $Z$ is generically finite over $q(Z)$.
\begin{proposition}\label{Prop: empty problematic locus implies geometrically non-degenerate}
    If $S_{\Lambda_Z}=\emptyset$, then $Z$ is geometrically non-degenerate.
\end{proposition}
\begin{proof}
    Suppose that $Z\subset A$ is geometrically degenerate. Thus, there is an abelian quotient $q:A\to X=A/Y$ such that $W\coloneqq q(Z)\subset X$ is a strict subvariety, and $\dim W<\dim Z$. We can assume that $W$ is geometrically non-degenerate. After lifting to an isogeny we can assume that $A=X\times Y$. Moreover,  we can assume that for all $x\in X$, 
    \[ \langle Z_x \rangle =Y \]
    (replace $Y$ with a smaller abelian subvariety and $X$ with a bigger one). We have
    \[ \Lambda_{W\times Y}\restr{Z}\subset \Lambda_Z \,. \]
    In particular, for all $x\in W$ and $\xi=N^\vee_{W,x}\in T^\vee X$, we have 
    \[ Y=\langle Z_x \rangle \subseteq \langle \Lambda_{Z,\xi} \rangle\,. \]
    As $W$ is geometrically non-degenerate, the Gauss map $\Lambda_W\to T^\vee_0 X$ is surjective \cite{Weissauer2011DegenerateGauss}. It follows that $T^\vee_0 X \subset S_{\Lambda_Z}$.
\end{proof}

The following example shows that the converse to Proposition \ref{Prop: empty problematic locus implies geometrically non-degenerate} is false (recall that divisors are geometrically non-degenerate if and only if they are ample \cite[188]{Debarre1995FultonHansen}):
\begin{example}\label{Ex: non-empty problematic locus}
   Define $\mathcal{A}^2_{t,g-t}$ to be the locus in $\mathcal{A}_g$ corresponding to ppav's $(A,\Theta)\in \mathcal{A}_g $ such that $A$ admits an abelian subvariety $X$ of dimension $t$ with $L_X\coloneqq \BO_A(\Theta)\restr{X}$ of type $(1,\dots,1,2)$ (see \cite{Debarre1988,AuffarthCodogni2019}). Let $2\leq t \leq g/2$ and let  $(A,\Theta)\in \mathcal{A}^2_{t,g-t}$ be general. Denote $Y$ the complementary abelian variety to $X$ (as in \cite[Sec. 5.3]{Birkenhake2004}). The polarization $L_Y\coloneqq \BO_A(\Theta)\restr{Y}$ is of type $(1,\dots,1,2)$ and the sum map induces a $4:1$ isogeny of polarized abelian varieties
    \[ \pi :X\times Y \to A \,. \] 
     Let $\tilde{\Theta}\coloneqq \pi^\ast \Theta \subset X\times Y$. We have $\tilde{\Theta}=\divv (s^X_1\otimes s^Y_1+s^X_2\otimes s^Y_2)$ where $s_1^X,s_2^X$ and $s_1^Y,s_2^Y$ are bases of $\HH^0(X,L_X)$ and $\HH^0(Y,L_Y)$ respectively. If $x\in B_X\coloneqq \divv( s^X_1)\cap \divv (s^X_2)$ and $y\in Y$ are general, then
     \[ N^\vee_{\tilde{\Theta},(x,y)}=\langle \dd_x (s_1^X)s^Y_1(y)+\dd_x(s_2^X)s^Y_2(y) \rangle \subset T^\vee_0 X \subset T^\vee_0 A\,. \]
     The Gauss map $\gamma_{\Lambda_{B_X}}:\Lambda_{B_X}\to T^\vee_0 X$ is surjective, so $\langle \Lambda_{\tilde{\Theta},\xi} \rangle=Y $ for a general $\xi\in  T^\vee_0 X$. Thus, $ T^\vee_0 X \subset S_{\Lambda_\Theta} $. By symmetry, we obtain
     \[ S_{\Lambda_\Theta}=  T^\vee_0 X \cup T^\vee_0 Y \,. \]
\end{example}

\subsection{Multiplicativity of the Chern-Mather class}\label{subsection: Mutliplicativity Chern-Mather class}
Let $p:A\times \PP V\to A$ denote the projection onto the first factor. For $\Lambda\in \mathscr{L}(A)$ and $r\geq0 $, the $r$-th \emph{Chern-Mather class} of $\Lambda$ is defined as
\begin{equation}\label{Equ: Definition of Chern-Mather class} 
c_{M,r}(\Lambda)\coloneqq p_\ast( h^{r}\cap [\PP \Lambda]) \in \HH_r(A) \,, 
\end{equation}
where $h$ is the pullback of the hyperplane class on $\PP V$. We define the total Chern-Mather class by
\[ c_{M}(\Lambda) \coloneqq \sum_{r=0}^g c_{M,r}(\Lambda) \in \HH_{\ast}(A)\,. \]
\begin{remark}
    The Chern-Mather classes could have been defined in the Chow group of $A$. We work with homology because we will want to apply $[n]_\ast$ to these classes and this works better on the level of homology with $\ZZ$-coefficients.
\end{remark}
We define the Chern-Mather class (resp. clean characteristic class) of a subvariety $Z\subset A$ by
\[ c_M(Z)\coloneqq c_M(\Lambda_Z) \,, \qquad \left(\text{resp.} \quad c_\ast(Z)\coloneqq c_M\circ \mathrm{cc}(\IC_Z) \right)\,.\]
Recall the following general properties of Chern-Mather classes of subvarieties of abelian varieties.
\begin{proposition}\label{Prop: Basic properties Chern Mather classes}
If $Z\subset A$ is a subvariety of dimension $d<g$, then
\begin{enumerate}
        \item $\deg c_\ast(Z)=\chi(\IC_Z)=\dim \omega_Z$.
        \item $c_{M,d}(\Lambda_Z)=[Z]$ and $c_{M,r}(\Lambda_Z)=0$ for $r\geq d+1$.
        \item Let $B\subset A$ be the largest abelian subvariety of $A$ such that $Z$ is stable under translation by $B$. Let $q:A\to A/B$ be the quotient, $d'\coloneqq \dim B$ and $W\coloneqq q(Z)$. We have
        \begin{align*}
            c_{M,k}(\Lambda_Z)&=0 &\text{for $0\leq k < d'$,}\\
            c_{M,k}(\Lambda_Z)&=q^\ast c_{M,k-d'}(\Lambda_{W})\neq 0 &\text{for $d'\leq k \leq d$.}
        \end{align*} 
        \item Let $D\subset A$ be a reduced, irreducible ample divisor, and let $d\coloneqq \dim \Sing(D)$. For $d<k\leq g-1$ we have
    \[ c_k(D)=[D]^{g-k}\,, \qquad \text{and} \qquad c_d(D)\neq [D]^{g-d} \,. \]
\end{enumerate}
\end{proposition}
\begin{proof}
1) follows from Kashiwara's Index Formula \cite[Th. 4.3.25]{Dimca2004} and (\ref{Equ: dim omega(P) is euler characteristic of P}). 2) follows from \cite[Lem. 3.1.2]{Kraemer2021MicrolocalGauss2}. For 3), let $B$ and $W$ be as in the proposition. We have $Z=q^{-1}(W)$, thus $c_{M,k}(\Lambda_Z)=q^\ast c_{M,k-d'}(\Lambda_{W})$ for all $k\geq d'$. Since $W$ is not stable under translation by a positive-dimensional abelian subvariety, the Gauss map $\gamma_{W}:\PP\Lambda_{W}\to \PP^{g-d'-1}$ is generically finite \cite{Weissauer2011DegenerateGauss}. Thus, by \cite[Lem. 3.1.2]{Kraemer2021MicrolocalGauss2} we have
\[ c_{M,k}(\Lambda_{W})\neq 0 \qquad \text{for $0\leq k \leq \dim W$}\,. \]
For 4), let $D$ be as in the proposition. We have
\[ \mathrm{CC}(\IC_D)=\Lambda_D+\sum_{Z\subset \Sing(D)} n_Z \Lambda_Z \,,\]
for some positive integers $n_Z\in \ZZ_{\geq0}$. By ii) the terms $\Lambda_Z$ do not contribute to $c_k(D)$ for $k>\dim \Sing(D)$ and contribute positively for $k\leq \dim \Sing(D)$. The Chern-Mather class of $\Lambda_D$ can be computed with Vogel's intersection algorithm, and the claim follows from \cite[Cor. 3.2.1]{Kraemer2021MicrolocalGauss2}. We obtain that there is an effective non-zero cycle $[V]$ with support the $d$-dimensional part of $\Sing(D)$, such that
\[ c_d(D)+[V]=[D]^{g-d}\,.\]

\end{proof}
Consider the quotient 
\[ \HH_{\leq d}(A)\coloneqq \faktor{\HH_{\ast}(A)}{\HH_{>d}(A)} \,, \quad \text{with} \quad \HH_{>d}(A) \coloneqq \bigoplus_{k> d} \HH_{k}(A)\,. \]
We define the Chern-Mather class in degree up to $d$ by
\[ c_{M,\leq d}(\Lambda) \coloneqq \sum_{k=0}^{d} c_{M,k}(\Lambda) \in \HH_{\leq d}(A) \,. \]
The Pontryagin product induces a ring structure on both $\HH_{\ast}(A)$ and $\HH_{\leq d}(A)$, which we denote by $\ast$. Let $\Lambda\in \mathfrak{L}(A)$. We would like to know when the Chern-Mather class is a ring morphism 
\[ c_M: (\langle \Lambda \rangle, \circ) \to (\HH_{\ast}(A),\ast) \qquad  \text{(resp. $c_M: (\langle \Lambda_Z \rangle, \circ) \to (\HH_{\leq d}(A),\ast)$).} \] 
The issue is that in our definition of the product $\circ$ on Lagrangian cycles, we exclude negligible components, whereas the Pontryagin product of the Chern-Mather classes computes the Chern-Mather class of the full non-clean product
\begin{equation}\label{Equ: Chern-Mather class of non-clean product is Pontryagin} c_{M}(\Lambda_1)\ast c_{M}(\Lambda_2)=c_{M}(\Lambda_1 \tilde{\circ} \Lambda_2)\,, \qquad \text{for $\Lambda_1,\Lambda_2\in \mathfrak{L}(A)$.} 
\end{equation}
Nevertheless, with our study of the problematic locus, we have a clear understanding of what can happen (compare with \cite[Lem. 3.3.1]{Kraemer2021MicrolocalGauss2}).
\begin{lemma}\label{Lem: Chern-Mather: multiplicativity relative to codim of problematic locus}
    Let $\Lambda \in \CL(A)$ and $d\coloneqq \codim_{ V} S_\Lambda$ (set $d=g+1$ if $S_\Lambda=\emptyset$). 
    \begin{enumerate}
        \item The Chern-Mather class is a ring morphism in degree up to $d-1$
    \[ c_{M,\leq d-1}:\left(\langle \Lambda \rangle,\circ \right) \to \left(\HH_{\leq d-1}(A) ,\ast \right)\,.\]
    
    \item If $\Lambda_1,\Lambda_2\in \langle \Lambda \rangle$ are effective cycles, then
    \[ c_M(\Lambda_1)\ast c_M(\Lambda_2)- c_M(\Lambda_1\circ \Lambda_2)\in \HH_{\geq d}(A) \]
    can be represented by an effective cycle.
    \item If $S_{\Lambda}$ is nonempty, then there exist $\Lambda_1,\Lambda_2\in \langle \Lambda \rangle$ such that the above cycle is nonzero.
    \end{enumerate}
\end{lemma}
\begin{proof}
i) follows from (\ref{Equ: negligible subvarieties}), (\ref{Equ: Chern-Mather class of non-clean product is Pontryagin}) and Proposition \ref{Prop: Basic properties Chern Mather classes}. ii) follows from (\ref{Equ: non-clean product is sum of clean and negligible with positive coefficients}). iii) follows from Lemma \ref{Lem: negligible components appear} and Proposition \ref{Prop: Basic properties Chern Mather classes}.
\end{proof}


\section{A Tannakian characterization of Jacobians}\label{A criterion for detecting Jacobians}
Recall that $\CL(A)$ is a $\lambda$-ring, and we denote by $\mathrm{Alt}^l:\CL(A)\to \CL(A)$ the $l$-th $\lambda$-operation. For $\Lambda\in \CL(A)$, we consider the generating polynomial
\[ E_\Lambda(x) \coloneqq \sum_{k\geq 0} c_{M}(\mathrm{Alt}^k(\Lambda)) x^k  \in \HH_\ast(A)\otimes_\ZZ \ZZ[x] \,, \]
The \emph{Eulerian polynomials} $P_n(x)\in \ZZ[x]$ are characterized by the formula
\[ \sum_{k\geq 1} k^n x^k = \frac{x P_n(x)}{(1-x)^{n+1}} \,, \quad \text{for $n\geq 0$\,.}\]
It is well-known that they satisfy the following inductive relation:
\begin{align*}
    P_0(x)&=1 \,, \\
    P_n(x)&=\sum_{k=0}^{n-1} \binom{n}{k} P_k(x)(1-x)^{n-1-k} \,, \quad  \text{for $n\geq 1$.}
\end{align*}
In particular, $P_n$ is a polynomial of degree $n-1$ for $n\geq 1$. We have the following analogue of \cite[Prop. 4.1]{Kraemer2016CharacteristicHPSeries} for the alternating product (where the product operation on $\HH_\ast(A)$ is the Pontryagin product).
\begin{lemma}\label{Lem: E lambda generating series Chern-Mather class}
    Let $\Lambda\in \CL(A)$, and assume that $c_{M,\leq d}:\langle \Lambda \rangle \to \HH_{\leq d}(A)$ is a ring morphism for some $d$. In $\HH_{\leq d}(A)\otimes_\ZZ \QQ[x]$ we then have
    \[ E_\Lambda(x)=(1+x)^{c_0} \prod_{k\geq 1} \exp\left(\frac{x P_{2k-1}(-x)}{(1+x)^{2k}} c_{M,k}(\Lambda) \right)\,, \]
    where $c_0=c_{M,0}(\Lambda)\in \ZZ $.
\end{lemma}
\begin{proof}
    Let $R=\QQ[e_1,e_2,\dots,e_n,\dots]$ be the ring of symmetric functions in infinitely many variables, where we denote by $e_n$ the $n$-th elementary symmetric function (see \cite[Sec. I.2]{macdonald1998symmetric}). For $n\geq 1$, let $p_n\in R$ be the $n$-th power sum. The generating series for the elementary symmetric functions can be expressed in terms of power sums by the following identity in $R[[x]]$:
    \[ \sum_{k\geq 0} e_k x^k = \exp\left( \sum_{k\geq 1} \frac{(-1)^{k+1}}{k} p_k x^k \right) \,, \]
    where $e_0=1$. By the theory of $\lambda$-rings (for instance \cite[Chap. I]{Knutson1973lambdarings}), we thus have in $\HH_{\leq d}(A)\otimes_\ZZ \QQ[x]$
    \begin{align*}
        E_\Lambda(x) &\coloneqq \sum_{k\geq 0} c_{M}(\mathrm{Alt}^k(\Lambda)) x^k \\
        &=\exp\left( \sum_{k\geq 1} \frac{(-1)^{k+1}}{k} c_{M}(\Psi^k(\Lambda)) x^k \right)\,, \\
    \end{align*}
    where $\Psi^k$ is the $k$-th Adams operation on $\mathscr{L}(A)$. By \cite[Prop. 1.3.3]{Kraemer2021MicrolocalGauss2}, we have $\Psi^k=[k]_\ast$ where $[k]:A\to A$ is the multiplication by $k$. Recall that $[k]_\ast$ is the multiplication by $k^d$ on $\HH_d(A,\QQ)$, so 

    \begin{align*}
        E_\Lambda(x) &= \exp\left( \sum_{k\geq 1} \sum_{l\geq 0} \frac{(-1)^{k+1}}{k} c_l k^{2l} x^k \right)\\
        &= (1+x)^{c_0} \prod_{l\geq 1} \exp\left(\frac{x P_{2l-1}(-x)}{(1+x)^{2l}} c_{l} \right)\,,
    \end{align*}
    where $c_l\coloneqq  c_{M,l}(\Lambda)\in \HH_{l}(A)$.
\end{proof}
In light of the above proposition, we make the following definition.
\begin{definition}\label{Def: coefficients Ekn i1 i2 etc}
    For $g,n\geq 0$, and $\underline{c}=(c_1,\dots,c_{g-1})$, we define
    \[ E^n(\underline{c},x)\coloneqq (1+x)^n \prod_{l=1}^{g-1} \exp\left(\frac{x P_{2l-1}(-x)}{(1+x)^{2l}} c_{l} \right)\in \QQ[[x]][c_1,\dots,c_{g-1}]\,.\]
    For $k\geq 0$ and $\underline{i}\in \ZZ_{\geq 0}^{g-1}$, we define the coefficients $E^n_k(\underline{i})$ by 
    \[ E^n(\underline{c},x)=\sum_{k\geq 0} \sum_{\underline{i}} E^n_k(\underline{i}) \underline{c}^{\underline{i}} x^k \,. \]
\end{definition}

Recall that we denote by $\bb^n_k\coloneqq \binom{n}{k}$ and $\bb_k\coloneqq \binom{2k}{k}$ the binomial and middle binomial coefficients, respectively. We make the following computation.
\begin{corollary}\label{Cor: first chern mather classes of Alt^n Lambda}
    Let $\Lambda\in \mathscr{L}(A)$ and $n\coloneqq \deg \Lambda$. If $c_{M,\leq 2}:\langle \Lambda \rangle \to \HH_{\leq 2}(A)$ is a ring morphism, then
\begin{align*}
    c_{M,0}(\mathrm{Alt}^k(\Lambda)) &= \bb^{n}_{k} \,, \\
    c_{M,1}(\mathrm{Alt}^k(\Lambda)) &= \bb^{n-2}_{k-1} c_1 \,, \\
    c_{M,2}(\mathrm{Alt}^k(\Lambda)) &= \bb^{n-4}_{k-2} \frac{c_1^{\ast 2}}{2}+\left(\bb^{n-4}_{k-3}-4\bb^{n-4}_{k-2}+\bb^{n-4}_{k-1}\right) c_2 \,,
\end{align*}
where $c_r\coloneqq c_{M,r}(\Lambda)$ for $r\geq 0$.
\end{corollary}
\begin{proof}
    In $\HH_{\leq 2}(A)\otimes_\ZZ \QQ[x]$ the formula of the lemma becomes
    \begin{align*}
        E_\Lambda(x)&= (1+x)^n \exp\left(c_1  \frac{x}{(1+x)^2}\right)  \exp\left(c_2  \frac{x(x^2-4x+1)}{(1+x)^4}\right) \\
        &= (1+x)^{n}+c_1 x (1+x)^{n-2}  + \left( \frac{c_1^{\ast2}}{2} x^2+ c_2 x(x^2-4x+1)\right)(1+x)^{n-4} \,,
    \end{align*}
    which gives the corollary after taking the coefficient of $x^k$.
\end{proof}
With the machinery introduced above, we now have an elegant way of computing the clean characteristic class of theta divisors of Jacobians:
\begin{proposition}\label{Prop: Chern-Mather classes of IC for Jacobians}
    Let $(JC,\Theta)\in \mathcal{J}_g$ be the Jacobian of a smooth curve $C$ of genus $g$. If $C$ is non-hyperelliptic, then
    \begin{align*}
        c_\ast(\Theta)&=\sum_{k=1}^{g} \bb_{k-1}[\Theta]^k/{k!} \cap[JC] \in \HH_\ast(JC)) \,.
        \end{align*}
    If $C$ is hyperelliptic, then
    \begin{align*} 
 c_\ast(\Theta) &= \sum_{k=1}^g \left( \bb_{k-1}-\bb^{2k-2}_{k-3}\right) [\Theta]^k\cap [JC] \in \HH_\ast(JC) \,.
    \end{align*}
\end{proposition}
\begin{proof}
We first treat the non-hyperelliptic case. Using the Abel-Jacobi map, we view $C$ as a subvariety of $JC$. By \cite[Th. 14]{Weissauer2006BrillNoetherSheaves} and \cite[124]{Weissauer2006BrillNoetherSheaves} (after choosing the right translate of $C$) there is an equivalence of Tannakian categories
\begin{align*}
   \omega: \langle \IC_C \rangle &\overset{\sim}{\longrightarrow} \mathrm{Rep}_\CC(\SL_{2g-2}(\CC)) \,,
\end{align*}
with $\omega(\IC_C)=[\CC^{2g-2}]$ and $\omega(\IC_\Theta)=[\mathrm{Alt}^{g-1}(\CC^{2g-2})]$. The Gauss map $\gamma_C:\PP \Lambda_C \to\PP^{g-1}$ is finite, so the Chern-Mather class is multiplicative (in all degrees) by Lemma \ref{Lem: Chern-Mather: multiplicativity relative to codim of problematic locus}. By Proposition \ref{Prop: Basic properties Chern Mather classes}, we have
\[ c_{0}(C)=2g-2\,, \qquad c_{1}(C)=[C] \eqcolon c\in \HH_1(JC) \,, \] 
and $c_{r}(C)=0$ for $r\geq 2$. We then have
\begin{align*}
    c_\ast(\Theta)&=c_M(\mathrm{Alt}^{g-1}(\Lambda_C))&\text{(\cite[Cor. 2.2.2]{Kraemer2021MicrolocalGauss2})} \\
    &= \left\lbrace (1+x)^{2g-2} \mathrm{exp}\left(x c/(1+x)^2\right) \right\rbrace_{x^{g-1}}  &\text{(Lemma \ref{Lem: E lambda generating series Chern-Mather class})} \\
    &=\sum_{k=0}^{g-1} \bb_{g-1-k} {c^{\ast k}}/{k!} & \\
    &= \sum_{k=0}^{g-1} \bb_{g-1-k}[\Theta]^k/(g-k)!\cap[JC] &\text{(\cite[322]{Birkenhake2004})}\,. 
\end{align*}
In the hyperelliptic case we have by \cite[p. 124]{Weissauer2006BrillNoetherSheaves} and \cite[Cor. 2.2.2]{Kraemer2021MicrolocalGauss2}
\[ \cc(\IC_\Theta)=\mathrm{Alt}^{g-1}(\Lambda_C)-\mathrm{Alt}^{g-3}(\Lambda_C) \,. \]
The result follows with a similar computation as in the non-hyperelliptic case.
\end{proof}

\subsection{Proof of Theorem \ref{maintheorem: Tannakian criterion for Jacobians}}\label{Subsec: proof of maintheorem Tannakian criterion}
   We can now prove Theorem \ref{maintheorem: Tannakian criterion for Jacobians}. Suppose that $(A,\Theta)\in \mathcal{A}_g$ is a fake Jacobian. By \cite[Th. 4.1.2]{Kraemer2021MicrolocalGauss2}, we have
    \[ [n]_\ast \cc(\IC_\Theta) = \mathrm{Alt}^{g-1}(\Lambda_Z) \qquad \text{(resp. $\mathrm{Alt}^{g-1}(\Lambda_Z)-\mathrm{Alt}^{g-3}(\Lambda_Z)$)}\,,\]

    for some $n> 0$ and $\Lambda_Z\in \langle \Lambda_\Theta \rangle$, in the non-hyperelliptic (resp. hyperelliptic) case. It is enough to show that $Z\coloneqq \mathrm{Supp}(\Lambda_Z)$ is of dimension $1$, since in that case $[n]_\ast \Theta$ (and thus $\Theta$) is a sum of curves and thus $(A,\Theta)$ is a Jacobian by \cite[Th. 1]{Schreieder2015ThetaCurveSummand}. We will do this by computing the Chern-Mather classes and apply Proposition \ref{Prop: Basic properties Chern Mather classes}. For $r\geq 0$, let $c_r\coloneqq c_{M,r}(\Lambda_Z)/n^{2r}\in \HH_{r}(A)\otimes_\ZZ \QQ$. By Lemma \ref{Lem: E lambda generating series Chern-Mather class}, there exist universal (graded) polynomials $Q_{g}(c_0,\dots,c_{g-1})$ (resp. $Q_{g}^{\mathrm{hyp}}$ in the hyperelliptic case) such that
    \begin{equation}\label{Equ: proof of theorem 1: relation chern class and universal polynomials}
    c_\ast(\Theta)=Q_g(c_0,\dots,c_{g-1}) \,, \qquad \left(\text{resp.} \quad  c_\ast(\Theta)=Q^{\mathrm{hyp}}_g(c_0,\dots,c_{g-1}) \,\right) 
    \end{equation}
    in the non-hyperelliptic (resp. hyperelliptic) case. We could compute these polynomials explicitly but we do not need to: by assumption, the left hand side of (\ref{Equ: proof of theorem 1: relation chern class and universal polynomials}) is the same as for Jacobians in degree up to $2$. We know that in the case of Jacobians, these relations hold with $c_2=\dots=c_{g-1}=0$. The equation in degree $0$ forces $c_0=2g-2$, and the equation in degree $1$ forces $c_1=[\Theta]^{\cdot(g-1)}/(g-1)!$. We just need to verify that the coefficients of $c_2$ in $Q_g$ and $Q_g^\mathrm{hyp}$ are non-zero. By Corollary \ref{Cor: first chern mather classes of Alt^n Lambda} these coefficients are
    \[\left\lbrace Q_g \right\rbrace_{c_2}= -\frac{2g+2}{g-2}\bb_{g-3} \,,\qquad  \left\lbrace Q_g^\mathrm{hyp} \right\rbrace_{c_2}=-4(2g-5)(g+3)\frac{(2g-6)!}{g!(g-3)!}\,. \]
    They are both non-zero, thus $c_2=0$ by (\ref{Equ: proof of theorem 1: relation chern class and universal polynomials}). Thus $Z$ is of dimension $1$, and we are done.

\section{The bielliptic Prym locus}\label{Sec: The bielliptic Prym locus}

\subsection{Notations}\label{Sec: Generalities fiber Gauss map on BEg}
We start by recalling general facts on the bielliptic Prym locus, following \cite{Debarre1988} and \cite{podelski2023GaussEgt}. Let $\pi:\tilde{C} \to C$ be an \'etale double cover, with $C$ a smooth genus $g+1$ bielliptic curve, i.e. admitting a cover $p:C\to E$ to an elliptic curve. We assume that the Galois group of $p\circ \pi$ is $(\ZZ/2\ZZ)^2$, thus inducing a tower of curves\footnote{The Pryms $\mathrm{Prym}(\tilde{C}/C)$ where the Galois group of $p\circ \pi$ is $\ZZ/4\ZZ$ are isomorphic to Pryms where the Galois group is $(\ZZ/2\ZZ)^2$ via the tetragonal construction \cite[Prop 3.8]{donagi1992fibers}.}
\begin{equation}\label{Fig: Cartesian Diagram of curves Egt} 
\begin{tikzcd}
   & \tilde{C} \arrow[dl,"\pi"'] \arrow[d,"\pi'"] \arrow[dr,"\pi''"] & \\
    C \arrow[dr,"p"'] & C' \arrow[d,"p'"] & C'' \arrow[dl,"p''"] \\
    & E & 
\end{tikzcd}\,.
\end{equation}
We can assume $g(C')=t+1 \leq g(C'')=g-t+1$ for some $0\leq t\leq g/2$. We define $\EE^\ast_{g,t}\subset \mathcal{A}_g$ as the set of Prym varieties $\mathrm{Prym}(\tilde{C}/C)$ obtained in this way and 
\[ \EE_{g,t}\coloneqq \overline{\EE_{g,t}^\ast}\subset \mathcal{A}_g\]
as the closure. The bielliptic Prym locus is the union
\[ \BE_g \coloneqq \bigcup_{t=0}^{\lfloor g/2 \rfloor} \EE_{g,t}\,. \] 
By \cite{Debarre1988}, the loci $\EE_{g,0},\dots,\EE_{g,\lfloor g/2 \rfloor}$ are distinct irreducible components of $\BE_g$. We define $\EE'_{g,t}\subset \EE_{g,t}$ to be the set obtained by relaxing slightly the conditions on $\EE_{g,t}^\ast$: $E$ is still assumed to be a smooth elliptic curve, but $\tilde{C}$ and $C$ are nodal curves, and $\pi$ is assumed to be admissible in the sense of Beauville \cite{Beauville1977}, i.e. the covering involution fixes only the nodes of $\tilde{C}$ and does not exchange the branches at the nodes. In the case of $\EE'_{g,t}$, the curves $C'$ and $C''$ are smooth (see the proof of \cite[Prop. 3.21]{podelski2023GaussEgt}). Moreover, we have $\EE'_{g,0}=\EE^\ast_{g,0}$. \par 
For $\ud=(d_1,\dots,d_n)$ a partition of $g$, we define $\SE_\ud\subset \EE_{g,0}$ to be the degenerations in $\EE_{g,0}$ where $\tilde{C}\to C$ is an admissible double cover of nodal curves and $E$ is a $n$-cycle of $\PP^1$'s, i.e. the $n$ irreducible components of $E$ are rational and the dual graph is the cyclic $n$-graph. For $g\geq 4$, we have \cite[Lem. 4.2.1]{podelski2023GaussEgt} 
 \begin{equation}\label{Equ: description bielliptic Prym locus} \EE_{g,t}=\left((\mathcal{J}_g \cup \mathcal{A}_g^\dec)\cap \EE_{g,t} \right)\cup \EE_{g,t}'\cup \bigcup_{ \ud \leq (t,g-t)} \SE_\ud \,,
\end{equation}
 where $\mathcal{A}_g^\dec\subset \mathcal{A}_g$ denotes the locus of decomposable ppav's and $\ud \leq \underline{e}$ if and only if $\ud$ is a subdivision of the partition $\underline{e}$. Let $(P,\Xi)\in \EE_{g,t}'$ or $(P,\Xi)\in \SE_\ud $ for some $\ud$. We use the notation of (\ref{Fig: Cartesian Diagram of curves Egt}). To the flat double cover $p: C\to E$ we can associate a line bundle $\delta\coloneqq \det( p_\ast \BO_C)^{-1} \in \Pic^g(E)$. If $\ddag$ stands for the symbol $\,',\,''$ or nothing, denote by $R^\ddag$ and $\Delta^\ddag$ the ramification and branch divisor of $p^\ddag$. We have $\Delta=\Delta'+\Delta''$, and $\Delta\in |\delta^{\otimes 2}|$. We define
 \[ \ns(\tilde{C}/C)=\# \{ (D',D'')\in E_t\times E_{g-t}\,|\, D'\leq \Delta'\,, D''\leq \Delta'' \,, \BO_E(D'+D'')= \delta \} \,,\]
 where we take $t=0$ if $(P,\Xi)\in \SE_\ud$. For $k\geq 0$, let
 \begin{align*}
     \EE_{g,t}^{\prime k}&\coloneqq \{ (P,\Xi)\in \EE_{g,t}'\,|\, \ns(\tilde{C}/C)=2k \} \,, \\
     \SE_\ud^k &\coloneqq \{ (P,\Xi)\in \SE_\ud'\,|\, \ns(\tilde{C}/C)=2k \} \,.
 \end{align*}
These loci are analogous to the $k$-th theta-null locus, because a general Prym in $\EE_{g,t}^{\prime k}$ or $\SE_\ud^k$ has exactly $k$ additional isolated singularities, which are at $2$-torsion points (see \cite{podelski2023GaussEgt}). \par 
 Let $(P,\Xi)\in \EE_{g,t}'$ for some $t\geq 0$. Let $\delta'\coloneqq \det(p'_\ast \BO_{C'})^{-1}\in \Pic(E)$, let $\Nm' :JC' \to JE$ be the norm map associated to $p'$, and
 \[ P' \coloneqq (\Nm')^{-1}(\delta') \subset \Pic(C^{\ddag}) \,. \]
 We define $\delta''$, $\Nm''$ and $P''$ similarly. By \cite[Prop. 5.5.1]{Debarre1988}, if $t=0$ (resp. $t\geq 1$), the pullback induces a $2:1$ (resp. $4:1$) isogeny
 \begin{equation}\label{Equ: 4:1 isogeny between P' P'' and P} 
 \pi'^\ast+\pi''^\ast:P'\times P'' \to P \,.
 \end{equation}

\subsection{The characteristic cycle of the Prym theta divisor}\label{Subsubsec: The characteristic cycle Egt smooth case}
In this section we compute the characteristic cycle $\mathrm{CC}(\IC_\Xi)$ for $(P,\Xi)\in \EE'_{g,t}$. Bressler and Brylinski, in \cite{BresslerBrylinski97}, show that if $(JC,\Theta)$ is the Jacobian of a smooth non-hyperelliptic curve, then the characteristic cycle of the intersection complex $\mathrm{CC}(\IC_\Theta)$ is irreducible. We will adapt the computation in \textit{loc. cit.} to our case.\par 
Let $g\geq 4$, $0\leq t\leq g/2$ and $(P,\Xi)\in \EE'_{g,t}$. We keep the notation of the previous section. Set-theoretically, we have
\begin{align*} 
P &= \{ L\in \Pic^{2g}(\tilde{C})\,|\, \Nm(L)=\omega_C\,,\, \hh^0(L) \text{ even} \} \\
 \Xi& = \{ L\in P \,|\, \hh^0(L)>0 \} \,. 
\end{align*} 
Let $\Theta'\subset \Pic^t(C)$, $\Theta''\subset \Pic^{g-t}(C'')$ be the theta divisors and $\alpha':C'_t\to \Theta'$, $\alpha'':C''_{g-t}\to \Theta''$ the Abel-Jacobi maps. Let
\[\Nm \coloneqq \Nm'+\Nm'':\Pic^t(C')\times \Pic^{g-t}(C'')\to \Pic^g(E)\,.  \] 
We define
\[ W\coloneqq( C'_t\times C''_{g-t}) \times_{\Pic^g(E)} \{\delta\} \,, \quad \text{and } \quad  \tilde{\Xi} \coloneqq (\Theta'\times \Theta'') \times_{\Pic^g(E)} \{\delta\} \,. \]
Let $h\coloneqq \alpha'\times \alpha''\restr{W}$, $g\coloneqq (\pi'^\ast+\pi''^\ast)\restr{\tilde{\Xi}}$ and $\phi\coloneqq g\circ h$. Then $g:\tilde{\Xi} \to \Xi$ is generically finite of degree $2$ \cite[Prop. 3.8]{podelski2023GaussEgt}. In summary, the following diagram is commutative:
\stepcounter{equation}
\begin{figure}[H]
\centering 
\begin{tikzcd}
\{\delta \} \arrow[d] & W \arrow[l] \arrow[r,"h"'] \arrow[d,phantom,sloped,"\subset"] \arrow[rr, bend left=20, "\phi"] & \tilde{\Xi} \arrow[r, "g"'] \arrow[d,phantom,sloped,"\subset"] & \Xi \arrow[r,hook] & \Pic^{2g}(\tilde{C})\\
{\Pic^g(E)}  & {C'_t\times C''_{g-t}} \arrow[l,"\Nm"] \arrow[r,"\alpha'\times\alpha''"'] & \Theta'\times \Theta'' \arrow[rr,hook] &  & \Pic^t(C')\times \Pic^{g-t}(C'') \arrow[u,"\pi'^\ast+ \pi''^\ast"']
\end{tikzcd}
\end{figure} 
\par 
For a curve $C$ the Brill-Noether loci are defined by 
\[W^r_k(C)\coloneqq \{L\in \Pic^k(C)\,|\, \hh^0(L)> r \}\,. \] 
We have the following:
\begin{proposition}\label{Prop: Fibers of g restr Xi tilde}
The locus where $g:\tilde{\Xi}\to \Xi$ fails to be finite is 
\[ \left((p'^\ast \Pic^1(E)+W^0_{t-2}(C'))\times (p''^\ast \Pic^1(E)+W^0_{g-t-2}(C'')\right)\times_{\Pic^g(E)}\{\delta\} \,. \]
On this locus, the fibers of $g$ are the orbits under action of the antidiagonal embedding $(p'^\ast,-p''^\ast):JE\to JC'\times JC''$.
\end{proposition}
\begin{proof}
In the proof of \cite[Prop 3.8]{podelski2023GaussEgt}, the only assumption made to prove the finiteness of $g$ at $(\BO_{C"}(D'),\BO_{C''}(D''))\in \tilde{\Xi}$ is that $D'$ or $D''$ is $p'$-simple (resp. $p''$-simple). From this, the proposition follows.
\end{proof}
Recall from \cite{BohroMacPherson1983PartResNilp} that a map of varieties $f:X\to Y$ is \emph{small} if for $k>0$,
\[ \codim \{y\in Y\,|\, \dim f^{-1}(y) \geq k \} <2k \,. \]
\begin{lemma}\label{Prop: phi : W to Xi is small}
Suppose $t\neq 2$ and $(P,\Xi)\in \EE'_{g,t}$, or $t=2$ and $(P,\Xi)\in \EE'_{g,2}$ general. Then the map $\phi:W\to \Xi$ is small.
\end{lemma}
\begin{proof}
Recall that since $C'$ is bielliptic, it is never hyperelliptic when $t>2$, and by our generality assumption it is not hyperelliptic if $t=2$. Thus
\begin{align*}
    W^1_{t}(C')&=\{p'^\ast \Pic^2(E)+W^0_{t-4}(C')\}\cup Z' \,, &\text{where $Z'\subset \Nm'^{-1}(\delta')$}\,, \\
    W^r_t(C')&=\{ p'^\ast \Pic^{r+1}(E)+W^0_{t-2r}(C')\}\,, & \text{for $r\geq 2$.}
\end{align*}
Similarly, we have
\[ W^1_{g-t}(C'')=\{p''^\ast \Pic^2(E)+W^0_{g-t-4}(C'')\}\cup Z'' \,, \quad \text{where $Z''\subset \Nm''^{-1}(\delta'')$.}\]
Let
\[ S'_k=p'^\ast \Pic^{k}(E)+W^0_{t-2k}(C')\,, \quad \text{and} \quad S''_k=p''^\ast \Pic^{k}(E)+W^0_{g-t-2k}(C'') \,, \]
for $k\geq 0$, where $S'_k=\emptyset$ (resp. $S''_k=\emptyset$) if $2k>t$ (resp. $2k>g-t$). Consider the following loci in $\Xi$
\begin{align*}
    S_n&\coloneqq g\left( \bigcup_{k} S'_{k}\times S''_{n-k}\cap \tilde{\Xi} \right) \,, \\
    T_n&\coloneqq g\left( (Z'\times S''_{n-1} \cup S'_{n-1}\cup Z'' )\cap \tilde{\Xi} \right) \,, \\
    Z&\coloneqq g(Z'\times Z'') \,.
    \end{align*}
For $k\geq 1$ have $\dim S_k=\dim T_k=g-2k$ and $\dim W\restr{S_k}=\dim W\restr{T_k}=g-k-1$. Thus, the general fiber of $\phi$ above $S_k\cup T_k$ is of dimension $k-1$. We have $\dim Z=g-6$ and $\dim W\restr{Z}=g-4$, so the general fiber above $Z$ is of dimension $2$. By \ref{Prop: Fibers of g restr Xi tilde}, away from $S_2\cup T_2\cup Z$, $\phi$ is birational, thus
\[ \{ x\in \Xi\,|\, \dim \phi^{-1}(x)\geq k \}= \begin{cases}
    S_{k+1}\cup T_{k+1}\cup Z & \text{if $k\in \{1,2\}$\,,}\\
    S_{k+1}\cup T_{k+1} & \text{if $3\leq k\leq g/2$.}
\end{cases}\]
In particular, we have for $1\leq k \leq g/2$
\[\codim  \{ x\in \Xi\,|\, \dim \phi^{-1}(x)\geq k \} >2k \,, \]
thus $\phi$ is small.
\end{proof}
By \cite[Prop. 3.9]{podelski2023GaussEgt}, we have
\[ \Sing(W)=\{ (D',D'')\in W\,|\, D'\leq R'\,, D''\leq R'' \} \,, \]
and these singularities are quadratic isolated singularities of maximal rank. Define $\,^t \dd \phi$ and $\phi_\pi$ to be the codifferential of $\phi$ and the projection onto the second factor respectively
\begin{equation}\label{Equ: definition of phi t and phi pi}
T^\vee W \overset{\,^t \dd \phi}{\longleftarrow} W\times_\Xi T^\vee P \overset{\phi_\pi}{\longrightarrow} T^\vee P \,. 
\end{equation}
Recall that for $x=(x',x'')\in P'\times P''$ there are canonical identifications $T^\vee_{g(x)} P =T^\vee_x(P'\times P'')$. By \cite[Sec. 3]{podelski2023GaussEgt} we have
\[ T^\vee_{x'} P' \simeq  \HH^0(E,\delta') \lhook\joinrel\xrightarrow{p'^\ast} \HH^0(C',\omega_{C'}))\simeq T^\vee_0 JC'\,,  \]
and similarly for $T^\vee_{x''}P''$. Under these identifications we have the following:
\begin{lemma}\label{Lemma: Kernel of codifferential of alpha in the smooth case}
If $D=(D',D'')\in W_\sm$, then
\[ \Ker(\,^t\dd_D \phi)=\left(\CC \cdot (s_{R'},s_{R''})+\bigoplus_{\ddag\in \{',''\}}\HH^0(C^\ddag,\omega_{C^\ddag}(-D^\ddag))\right) \cap (\oplus_\ddag \HH^0(E,\delta^\ddag))\,, \]
where $s_{R^\ddag}\in \HH^0(C^\ddag,\omega_{C^\ddag})$ corresponds to the pullback to $C^\ddag$ of the same generator $\dd z$ of $\HH^0(E,\omega_E)$.
\end{lemma}
\begin{remark}
    Note that under the identification 
    \[ T^\vee JC'\oplus T^\vee JC''=(T^\vee JE \oplus T^\vee P')\oplus(T^\vee JE \oplus T^\vee P'') \,, \]
    $\CC \cdot ( s_{R'},s_{R''})$ corresponds to the diagonal embedding of $T^\vee JE \hookrightarrow T^\vee JE \oplus T^\vee JE$.
\end{remark}
\begin{proof}
From the definition of $\phi$ we derive the following maps of vector bundles on $W$ (we omit the pullback to $W$ in the notation)
\[ \begin{tikzcd}
     T^\vee W & T^\vee C'_t\times T^\vee C''_{g-t} \arrow[l,two heads,"{^t\dd \iota_W}"] & T^\vee JC'\oplus T^\vee JC''\arrow[l,"{^t\dd \alpha}"] & T^\vee P' \oplus T^\vee P'' \arrow[l,hook'] \arrow[lll,"^t\dd \phi"',bend right=15]
\end{tikzcd}
\]
It is well known (see \cite[Lem 2.3 p. 171]{arbarello}) that the kernel of $^t \dd \alpha$ at $D=(D',D'')\in W$ is identified with
\[ \Ker( ^t \dd\alpha)_D= \HH^0(C',\omega_{C'}(-D'))\oplus \HH^0(C'',\omega_{C''}(-D'')) \,. \]
Under the assumption that $D\in W_\sm$ we have $D\nleq(R',R'')$ and 
\[ 0\neq \Ker(\,^t\dd \iota_W) = \langle \,^t \dd \alpha(s_R',s_R'') \rangle \,. \]
The lemma follows.
\end{proof}

\begin{lemma}\label{Lem: pushforward conomral W to Xi}
Let $(P,\Xi)\in \EE'_{g,t}$ and 
\[ \Lambda \coloneqq \,^t \dd \phi^{-1}(\Lambda_{W_\sm}) \subset W_\sm \times_P T^\vee P \]
where $\Lambda_{W_\sm}$ is the zero section of $N^\vee W_\sm$. If $t=0$, $\Lambda$ is irreducible and of degree $2$ above $\Lambda_\Xi$. If $0<t\leq g/2$, then $\Lambda=\Lambda' + \Lambda_{W'\times W''}$, where $\Lambda'$ is of degree $2$ above $\Lambda_\Xi$.
\end{lemma}
\begin{proof}
Let $D=(D',D'')\in W_\sm$. We split the proof into four cases.
\begin{enumerate}
    \item $D'$, $D''$ are $p^\ddag$-simple and $\Nm(D^\ddag)\neq \delta^\ddag$ for $\ddag\in \{',''\}$. Then by \cite[Prop. 3.3]{podelski2023GaussEgt} we have 
    $\hh^0(C^\ddag,D^\ddag)=\hh^0(C^\ddag,\omega_{C^\ddag}(-D^\ddag))=1$ and the corresponding non-zero section can be written as
    \[ s^\ddag=a^\ddag s_{R^\ddag}+b^\ddag p^{\ddag,\ast}t^\ddag \]
    where $a^\ddag,b^\ddag\neq 0$ and $t^\ddag\in \HH^0(E,\delta^\ddag)$. Thus, by Lemma \ref{Lemma: Kernel of codifferential of alpha in the smooth case}
    \[\Ker(\,^t \dd\phi)_D=\langle( \frac{a'}{b'} t',\frac{a''}{b''}t'')\rangle\in \HH^0(E,\delta')\oplus\HH^0(E,\delta'') \,. \]
    Points in this case correspond to an irreducible open set of $\Lambda$, which is clearly of degree $2$ over $\Lambda_\Xi$. This locus is of dimension $g$.
    
    \item $D'$, $D''$ are $p^\ddag$-simple and $\Nm(D')=\delta'$ (then necessarily $\Nm(D'')=\delta''$). If $\hh^0(C',D')=1$ (resp. $\hh^0(C'',D'')=1$), then using the above notation we have $a'=0$ (resp. $a''=0$) and
    \[\Ker(\,^t \dd\phi)_D=\langle (t',0),(0,t'')\rangle\in \HH^0(E,\delta')\oplus\HH^0(E,\delta'') \,, \]
    where $\divv t^\ddag=p^\ddag_\ast D^\ddag$. Let $\Lambda_{W'\times W''}$ be the closure of points which fall in this case. Clearly this locus is empty if $t=0$. $\Lambda_{W'\times W''}$ maps to $W'\times W'' \subset W$, where $W^\ddag=\Nm^{\ddag,-1}(\delta^\ddag)\subset C^\ddag_{g^\ddag-1}$ which is of dimension $(t-1)+(g-t-1)$ and the generic fiber is of dimension $2$, so $\Lambda_{W'\times W''}$ is of dimension $g$.
    
    \item $D'$, $D''$ are $p^\ddag$-simple, $\Nm(D^\ddag)=\delta^\ddag$ and $\hh^0(C^\ddag,D^\ddag)>1$. Then $\hh^0(C^\ddag,D^\ddag)=2$ by \cite[Prop. 3.3]{podelski2023GaussEgt}, and by the above reasoning $\dim \Ker(\,^t \dd\phi)_D=3$. The set of points obtained in this way projects to a locus of dimension $(t+1-4)+(g-t+1-4)$ in $W$ and the fibers are of dimension $3$, thus the total dimension of this set is $g-3$.\\
    
    \item Finally, suppose $D^\ddag=p^{\ddag,\ast} M^\ddag+F^\ddag$ where $F^\ddag$ is $p^\ddag$-simple, $\deg M^\ddag=k^\ddag$ and $k'+k''>0$. Say for example $k'>0$. Then
    \[ \HH^0(C',\omega_{C'}(-D'))=\HH^0(E,\delta'(-M'-p'_\ast D'))\subset T^\vee P'\,. \]
    Thus by Lemma \ref{Lemma: Kernel of codifferential of alpha in the smooth case} we have
    \[ \Ker(\,^t \dd \phi)_D = \bigoplus_{\ddag\in \{',''\}} \HH^0(E,\delta^\ddag(-M^\ddag-p^\ddag_\ast D^\ddag)) \,. \]
    Let $W_k$ be the set of points obtained in this way with $k'+k''=k$. Then $W_k$ projects to a set of dimension $k'+(t-2k')+k''+(g-t-2k'')-1$ in $W$ and the fibers of the projection are of dimension $k'+k''$ thus $W_k$ is of dimension $g-1$.
\end{enumerate}
Since every irreducible component of $\Lambda$ is of dimension at least $g$, the result follows.
\end{proof}
We can now prove the following:
\begin{theorem}\label{Theorem: Characteristic Cycle for E'gt}
Let $g\geq 4$, $0\leq t\leq g/2$, $(P,\Xi)\in \EE'_{g,t}$ and suppose that $\Nm(D')\neq \delta'$ for all $(D',D'')\in W_\sing$. Then the characteristic cycle of the intersection complex of $\Xi$ is given by
\[ \mathrm{CC}(\IC_\Xi)=
\begin{cases}
    \Lambda_\Xi & \text{if $g$ is even,}\\
    \Lambda_\Xi + \sum_{x\in \Xiasing}  \Lambda_x & \text{if $g$ is odd,}
\end{cases} \]
where $\Xiasing=\phi(W_\sing)$ is the set of additional isolated singularities.
\end{theorem}
\begin{proof}
Let $\Xi^\mathrm{o}=\Xi\setminus \Xiasing$. By Proposition \ref{Prop: phi : W to Xi is small}, the map $\phi^\mathrm{o}: W_\sm\to \Xi^\mathrm{o}$ is small. Thus, by the BBDG decomposition Theorem,
\[ R\phi^\mathrm{o}_{\ast}(\underline{\CC}_{W_\sm})=\IC_{\Xi^\mathrm{o}} \oplus \mathscr{F} \,, \]
for some perverse sheaf $\mathscr{F}$ on $\Xi^\mathrm{o}$. Thus, by \cite[Prop. 5.4.4]{KashiwaraShapira1990} we have an inclusion
\[ \mathrm{CC}(\IC_{\Xi^\mathrm{o}}) \subset \phi^\mathrm{o}_{\pi,\ast}(\,^t \dd\phi^{-1} (\Lambda_{W_\sm})) \,. \]
Thus, by Lemma \ref{Lem: pushforward conomral W to Xi} we have
\[ \mathrm{CC}(\IC_{\Xi})=\Lambda_{\Xi}+m\phi_\pi(\Lambda_{W'\times W''})+\sum_{x\in \Xiasing} m_x \Lambda_x \]
for some integers $m,m_x$. By \cite[Prop. 5.2.2]{Debarre1988}, a general point in $\phi(W'\times W'')$ is a smooth point of $\Xi$, so $m=0$. Finally, by \cite[Prop 3.10]{podelski2023GaussEgt}, under the assumption of the theorem, the points in $\Xiasing$ are isolated quadratic singularities of maximal rank of $\Xi$. It is well known that in this case $m_x=0$ if $g$ is even and $m_x=1$ if $g$ is odd \cite[Lem. 6.5.2]{Kraemer2021MicrolocalGauss2}.
\end{proof}

\subsection{The Gauss map for bielliptic Pryms }
We now describe the fibers of the Gauss map more precisely in the case $t=0$. 

\begin{lemma}\label{Lem: Gauss map finite for Eg0'}
    Let $g\geq 4$ and $(P,\Xi)\in \EE^{\prime 0}_{g,0}$. The Gauss map $\gamma_\Xi: \PP \Lambda_\Xi \to  \PP T^\vee_0 P$ is finite.
\end{lemma}
\begin{proof}
Let $(P,\Xi)=\Prym(\tilde{C}/C)\in \EE'^0_{g,0}$. We keep the notation of the two previous sections.  Recall that in this setting, all the curves of (\ref{Fig: Cartesian Diagram of curves Egt}) are smooth. We have canonical identifications
\[ \HH^0(E,\delta'')\simeq T^\vee_0 P'' \simeq T_0^\vee P \,. \] 
The pullback of divisors induces an inclusion $p''^\ast:|\delta''|\hookrightarrow |\omega_{C''}|$. Let $\mathcal{F}:|\omega_{C''}|\to |\delta''|$ be the projection from $R''\in|\omega_{C''}|$. Let $\GG'':\Theta''\dashrightarrow |\omega_{C''}|$ and $\GG:\Xi \dashrightarrow |\delta''|$ be the respective Gauss maps. The following diagram is commutative \cite[Fig. 2.16]{podelski2023GaussEgt}:
\begin{equation*}
\begin{tikzcd}
W \arrow[r] \arrow[rrr,dashed,bend left=15,"\tilde{\GG}"] \arrow[d,phantom,sloped,"\subset"] & \Xi'' \arrow[r, "g"'] \arrow[d,phantom,sloped,"\subset"] & \Xi \arrow[r,dashed, "\GG"'] & {|\delta''|} \\
C''_{g}  \arrow[rrr,bend right=15,dashed," \tilde{\GG}''"']\arrow[r,"\alpha''"] &  \Theta''\arrow[rr,dashed," \GG''"] & & { |\omega_{C''}|} \arrow[u,dashed,"\F"']
\end{tikzcd}
\end{equation*}

For $D\in C''_g$ with $\hh^0(C'',D)=1$, $\tilde{\GG}''(D)$ is the unique canonical divisor $H\in |\omega_{C''}|$ such that $D\leq H$ \cite[p. 247]{arbarello}. For $M\in |\delta''|$, let 
\[ V_M\coloneqq  \overline{\mathcal{F}^{-1}(M)} \subset |\omega_{C''}| \,. \] 
$V_M$ is the line (in the projective space $|\omega_{C''}|\simeq \PP^g$) generated by $R'',p''^\ast M\in |\omega_{C''}|$. Let
\[ \Lambda \coloneqq \{(D,M)\in W\times |\delta''|\,\quad \text{with}\quad  D\leq H \quad \text{for some $H\in V_M$} \} \,. \] 
Let 
\[ Z_M\coloneqq \overline{\GG''^{-1}(V_M)}=\{ D\in C''_g\,|\, \exists H\in V_M\,, D\leq H \} \,. \]
Clearly $Z_M\to V_M$ is finite. Recall that 
 \[ W\coloneqq (\Nm''\circ \alpha'')^{-1}(\delta) \subset C''_g \,. \] 
 By the definition of $\mathscr{E}'^0_{g,0}$, we have $\ns(\tilde{C}/C)=0$, so
 \[ \{D \in W \,|\, D\leq R'' \}=\emptyset \,. \] 
 Thus, no point of $Z_M \cap W$ lies above $R''$. Thus, $Z_M\cap W=\Lambda\restr{M}$ is finite. We have $\PP\Lambda_\Xi = g\times \Id_{|\delta''|}(\Lambda)$ and thus $\gamma_\Xi$ is finite.
\end{proof}

 \subsection{The Chern-Mather class of the Prym theta divisor}
Theorem $1$ of \cite{podelski2023GaussEgt} corresponds to the computation of the $0$-th Chern-Mather class $c_{M,0}(\Lambda_\Xi)$. With the same arguments we can also compute the higher-dimensional Chern-Mather classes.
\begin{lemma}\label{Lem: Chern Mather classes for E'g0}
    If $g\geq 4$ and $(P,\Xi)\in\EE'_{g,0}$, then for $r\geq 1$ 
    \[ c_{M,r}(\Lambda_\Xi)=\binom{2g-2r-2}{g-r-1}\frac{ [\Xi]^{g-r}}{(g-r)!} \in \HH_{r}(P)\,. \]
\end{lemma}
\begin{proof}
    We keep the notation of the proof of Lemma \ref{Lem: Gauss map finite for Eg0'}. The proof of Theorem $1$ in \cite{podelski2023GaussEgt} starts by taking a general point $M\in |\delta''|$ and computing its preimage by the Gauss map. To compute the Chern-Mather classes, we replace the point $M$ with a general $r$-plane $H\subset |\delta''|$. Let
    \begin{align*} V_H&= \overline{\FF^{-1}(H)} = \langle R'' , p''^\ast H \rangle \subset |\omega_{C''}| \,, \\
    Z_H&=\{D\in C''_g \,|\, D\leq F\in V_H \}\subseteq C''_g\,. 
    \end{align*}
    $Z_H$ is the degree $g$ secant variety of the linear series $V_H$ and thus by \cite[p. 342]{arbarello} we have
\[ [Z_H]= \sum_{k=0}^{g-r-2} \binom{g-r-2}{k} \frac{ x^k \theta^{g-r-1-k}}{(g-r-1-k)!}\in \HH_{r+1}(C''_g) \,, \]
where $x=[C''_{g-1}]\in \HH^1(C'')$ and $\theta=[\Theta'']\in \HH^1(JC'')$. The locus of indeterminancy of $\FF$ is $\{R''\}$ and thus, for dimension reasons, when $r\geq 1$ no component of $Z_H\cap W$ is in the locus of indeterminancy of the Gauss map $\tilde{\GG}: W\dashrightarrow |\delta''|$. Thus
\begin{align*}
    c_{M,r}(\Xi'')&=\left[\overline{\GG''^{-1}(V_H)}\right]\\
    &=\alpha''_\ast [Z_H]\cdot [W] \\
    &= \frac{\theta^{g-r}}{(g-r)!}\binom{2g-2r-2}{g-r-1} \cap [P''] \\
    &=\frac{[\Xi'']^{g-r}}{(g-r)!} \binom{2g-2r-2}{g-r-1} \in \HH_{r}(P'')\,.
\end{align*}
We conclude by the fact that $\pi''^\ast\restr{P''}:P''\to P$ is an isogeny of polarized abelian varieties and $(\pi''^\ast\restr{P''})^\ast \Xi=\Xi''$.
\end{proof}
\begin{remark}
    Another proof of the above result is the following. If $(P,\Xi)\in \EE'_{g,0}$, then $\IC_{\Theta''}$ is non-characteristic with respect to the embedding $i:P''\hookrightarrow JC''$ (in the sense of \cite[Def. 5.4.12]{KashiwaraShapira1990}) if and only if $\ns(\tilde{C}/C)=0$. In that case $\IC_{\Xi''}=i^\ast \IC_{\Theta''}[-1]$ and by \cite[Prop. 9.4.3, Cor. 10.3.6]{KashiwaraShapira1990} we have
    \[\mathrm{CC}(\IC_{\Xi''})=\,^t \dd i( i_\pi^{-1}(\mathrm{CC}(\IC_{\Theta''})))=\Lambda_{\Xi''}\,, \] 
     where $i_\pi: T^\vee JC''\restr{P''}\hookrightarrow T^\vee JC''$. It follows that
    \[(\pi''^\ast\restr{P''})^\ast\left(c_M(\mathrm{cc}(\IC_\Xi))\right)=c_M(\cc(\IC_{\Xi''}))=c_M(\cc(\IC_{\Theta''}))\cap [P''] \,, \]
    so Lemma \ref{Lem: Chern Mather classes for E'g0} follows from Proposition \ref{Prop: Chern-Mather classes of IC for Jacobians}. This explains why $c_M(\Lambda_\Xi)$ ``looks like'' the Chern-Mather class of a Jacobian.
\end{remark}

We now do a similar computation for the case of $(P,\Xi)\in\EE_{g,t}'$ with $t\geq 1$. For $\ddag\in \{',''\}$, we have by \ref{Equ: 4:1 isogeny between P' P'' and P} an inclusion $\pi^{\ddag\ast}:P^\ddag\hookrightarrow P$ and a $4:1$ isogeny $P'\times P''\to P$. Let $\Nm_{P^\ddag}:P\to P^\ddag$ denote the corresponding norm map. Let $\xi=[\Xi]\in \HH^1(P)$ denote the principal polarization on $P$. Let
\begin{align*}
    \xi'&\coloneqq [\Theta']\cap P'\in \HH^1(P')\,, & \xi''&\coloneqq [\Theta'']\cap P''\in \HH^1(P'') \,, \\
   \tilde{\xi}'&\coloneqq \frac{1}{4}\Nm_{P'}^\ast \Xi' \in \HH^1(P)_\QQ\,, & \tilde{\xi}'' &\coloneqq \frac{1}{4}\Nm_{P''}^\ast \Xi'' \in \HH^1(P)_\QQ\,.
\end{align*}
The composite $(\Nm_{P'}, \Nm_{P''})\circ (\pi'^\ast + \pi''^\ast)$ is the multiplication by $2$ on $P'\times P''$ thus
\[ (\pi'^\ast + \pi''^\ast)^\ast(\tilde{\xi}'+\tilde{\xi}'')=\xi+\xi'=(\pi'^\ast + \pi''^\ast)^\ast \xi\,.\]
Thus
\begin{equation*}
    \xi= \tilde{\xi}'+\tilde{\xi}'' \,, \quad  \pi'\restr{P'}^\ast \tilde{\xi}'=\xi' \,, \quad \pi''\restr{P''}^\ast  \tilde{\xi}''=\xi'' \,. 
\end{equation*}
With the above notations, we have the following:
\begin{lemma}\label{Lem: Chern Mather class of E'gt t>=1}
    Let $t\geq 1$ and $(P,\Xi)\in \EE'_{g,t}$. Then for $r\geq 1$, the $r$-th Chern-Mather class of $\Xi$ is given by
    \[ c_{M,r}(\Lambda_\Xi)=\sum_{k=0}^{g-r} \left(\bb_k \bb_{g-r-k-1}+\bb_{k-1}\bb_{g-r-k} \right) \frac{\tilde{\xi}^{\prime k}}{k!} \cdot \frac{\tilde{\xi}^{\prime \prime g-r-k}}{(g-r-k)!}\cap [P] \,, \]
   where $\bb_k\coloneqq \binom{2k}{k}$ denotes the middle binomial coefficient.
\end{lemma}
\begin{proof}
     We have a similar commutative diagram as in the proof of Lemma \ref{Lem: Gauss map finite for Eg0'} (see \cite[Fig. 3.12]{podelski2023GaussEgt})
     
\begin{equation*}
\begin{tikzcd}
\{\delta\} \arrow[d,hook] & W  \arrow[l,hook'] \arrow[r]  \arrow[rrr,dashed,bend left=15,"\tilde{\GG}"] \arrow[d,phantom,sloped,"\subset"] & \tilde{\Xi} \arrow[r, "g\restr{\tilde{\Xi}}"'] \arrow[d,phantom,sloped,"\subset"] & \Xi \arrow[r,dashed, "\GG"'] & \PP^{g-1} \\
\Pic^g(E) & C'_t\times C''_{g-t} \arrow[l,"\Nm",two heads] \arrow[rrr,bend right=15,dashed,"\tilde{\GG'}\times \tilde{\GG''}"']\arrow[r,"\alpha"] & \Theta'\times \Theta''\arrow[rr,dashed,"\GG'\times \GG''"] & & {|\omega_{C'}|\times |\omega_{C''}|} \arrow[u,dashed,"\F"']
\end{tikzcd}
\end{equation*}

Let $r\geq 1$ and $H\subset \PP^{g-1}$ be a general $r$-plane. Let $Z\coloneqq \overline{C'_t\times C''_{g-t}\restr{H}}$. Because $\dim Z\geq 2$, no component of $ W \cap Z$ is contained in the locus where $\tilde{\GG}$ is undefined. The morphism $g\restr{\tilde{\Xi}}$ is of degree $2$ \cite[Prop. 3.8]{podelski2023GaussEgt}, and thus by the above diagram we have
\[ c_{M,r}(\Lambda_\Xi)=\frac{1}{2} g_\ast \alpha_\ast( [W]\cdot [Z])\in \HH_{r}(P) \,. \] 
By \cite[Equ. 3.15]{podelski2023GaussEgt} we have
\[ \alpha_\ast [Z]= \sum_{k=0}^{g-r-1} \bb_k \bb_{g-r-k-1} \frac{\theta^{\prime k+1}}{(k+1)!}\cdot \frac{\theta^{\prime \prime g-r-k}}{(g-r-k)!} \cap [JC'\times JC'']\,, \]
where $\theta'$ and $\theta''$ correspond to the principal polarization on $JC'$ and $JC''$ respectively.
Similarly to how we defined $\tilde{\xi}'$ and $\tilde{\xi}''$, let
\begin{align*}
    {\theta}'_E&\coloneqq \frac{1}{2}\Nm'^{\ast}(\theta_E)\,, &  {\theta}''_E&\coloneqq \frac{1}{2}\Nm''^{\ast}(\theta_E)\,, \\
    \check{\xi}'& \coloneqq \frac{1}{4}\Nm_{JC'\to P'}^\ast \xi'=(\pi'^\ast)^\ast \tilde{\xi}' \,, & \check{\xi}''& \coloneqq \frac{1}{4}\Nm_{JC''\to P''}^\ast \xi'' =(\pi''^\ast)^\ast \tilde{\xi}''\,, \\
\end{align*}
where $\theta_E$ denotes the principal polarization on $E$. With the same argument, we have
\[ \theta'={\theta}'_E+\check{\xi}'\,, \quad \theta''={\theta}''_E+\check{\xi}''\,.\]
We have $\Nm^\ast \theta_E= 2{\theta}'_E+2{\theta}''_E+x$ where $x\in \HH^1(JC',\QQ)\otimes \HH^1(JC'',\QQ)$ under the Künneth decomposition. Any class in $\HH_\ast(W,\QQ)$ which does not have ${\theta}'_E{\theta}''_E$ as a factor is annihilated by $\pi'^\ast+\pi''^\ast$, and therefore
\begin{align*}
    c_{M,r}&(\Lambda_\Xi)=\frac{1}{2} g_\ast \alpha_\ast ([W]\cdot [Z])\\
    &= \frac{1}{2}g_\ast \Big(2({\theta}'_E+{\theta}''_E) \sum_{k} \frac{\bb_k \bb_{g-r-k-1}}{k!(g-r-k-1)!}\big((k+1){\theta}'_E \check{\xi}'^k+\check{\xi}'^{k+1} \big)  \\
        & \qquad \cdot  \big((g-r-k){\theta}''_E \check{\xi}'^{g-r-k-1}+\check{\xi}'^{g-r-k} \big) \Big)\cap [JC'\times JC''] \\
        &=g_\ast \sum_k \frac{ \bb_k \bb_{g-r-k-1}}{4}\Big( \frac{\check{\xi}'^k \check{\xi}''^{g-r-k}}{k!(g-r-k)!} +\frac{\check{\xi}'^{k+1}\check{\xi}''^{g-r-k-1}}{(k+1)!(g-r-k-1)!} \Big)\cap [P'\times P'']\\
        &=\sum_{k=0}^{g-r} \left(\bb_k \bb_{g-r-k-1}+\bb_{k-1}\bb_{g-r-k} \right) \frac{\tilde{\xi}^{\prime k}}{k!} \cdot \frac{\tilde{\xi}^{\prime \prime g-r-k}}{(g-r-k)!}\cap [P]\,. 
\end{align*}

\end{proof}

 \section{Proof of Theorem \ref{maintheorem: Tannakian Schottky on Bielliptic Prym locus}}\label{Sec: The proof of Theorem 3}
 Recall that by Proposition \ref{Prop: Basic properties Chern Mather classes}, the dimension of the Tannakian representation is given by
 \[ \dim \omega_\Theta= \deg c_\ast(\Theta)=\deg(\cc(\IC_\Theta)) \,, \qquad \text{for $(A,\Theta)\in \mathcal{A}_g$}\,. \] 
The first step in the proof of Theorem \ref{maintheorem: Tannakian Schottky on Bielliptic Prym locus} is to find those Pryms whose characteristic class has the right degree.
\begin{lemma}\label{Lem: which bielliptic Pryms have same Euler Char as Jac}
    If $g\geq 4$ and $(P,\Xi) \in \BE_g \setminus (\mathcal{J}_g\cap \BE_g)$, then
    \[ \dim \omega_\Xi =\bb_{g-1} \iff (P,\Xi)\in \EE^{\prime 0}_{g,0} \cup \SE_{g}^0 \cup \SE_{1,g-1}^0 \,. \]
\end{lemma}
\begin{proof}
    Let $g\geq 4$. Suppose first $(P,\Xi)\in \EE_{g,t}$ general with $2\leq t \leq g/2$. By Theorem \ref{Theorem: Characteristic Cycle for E'gt}, we have
    \[ \mathrm{CC}(\IC_\Xi)=\cc(\IC_\Xi)=\Lambda_\Xi \,.\]
    By \cite[Th. 1 and Prop 3.21]{podelski2023GaussEgt}, we have
    \begin{align}\label{Equ: Euler characteristic IC Xi for Egt}
    \deg c_\ast(\Xi)=\deg \Lambda_\Xi = \bb_{t-1}\bb_{g-t}+\bb_{t}\bb_{g-t-1}-2^{g-1} < \bb_{g-1}\,.
    \end{align}
    By \cite[Th. 1.6]{KraemerCodogni}, $\deg c_\ast(\Xi)$ only decreases under specialization (the argument made for the degree of the Gauss map extends to $\deg c_\ast(\Theta)$). Thus, inequality (\ref{Equ: Euler characteristic IC Xi for Egt}) holds on all $\EE_{g,t}$ for $t\geq 2$. Let $\ud$ be a partition of $g$ distinct from $(g)$ and $(1,g-1)$. Then $\ud\leq (t,g-t)$ for some $2\leq t \leq g/2$ and by (\ref{Equ: description bielliptic Prym locus}), we have $\SE_\ud\subset \EE_{g,t}$ thus 
    \[ \deg c_\ast(\Xi) <\bb_{g-1} \,, \qquad \text{for all $(P,\Xi)\in \SE_\ud$.} \]
Let $k\geq 0$ and $(P,\Xi)\in \EE_{g,0}^{\prime k}\cup \SE_g^k\cup \SE_{1,g-1}^k$. Then, by Theorem \ref{Theorem: Characteristic Cycle for E'gt} and \cite[Th. 3]{podelski2023boundary} we have
\[ \mathrm{CC}(\IC_\Xi)=\begin{cases}
    \Lambda_\Xi & \text{if $g$ even,}\\
    \Lambda_\Xi + \sum_{x\in \Xiasing} \Lambda_x &\text{if $g$ odd,}
\end{cases}\]
where $\#(\Xiasing)=k$. Thus, by \cite[Th. 1 and 2]{podelski2023GaussEgt} we have
\[ \deg c_\ast(\Xi)= \begin{cases}
    \bb_{g-1}-2k &\text{if $g$ even,} \\
    \bb_{g-1}-k & \text{if $g$ odd.}
\end{cases}\]
Thus, $\deg c_\ast(\Xi)=\bb_{g-1}$ if and only if $k=0$. \par 
It remains to treat the case of $\EE'_{g,1}$. If $g\geq 5$ and $(P,\Xi)\in \EE'_{g,1}$, then by \cite[Cor. 5.7]{podelski2023GaussEgt}
    \[ \deg c_\ast(\Xi) \geq \deg \Lambda_\Xi > \bb_{g-1}\,. \]
We now deal with the case $g=4$. For $(P,\Xi)\in \mathcal{A}_4$ we have by \cite[Sec. 6.4]{Kraemer2021MicrolocalGauss2}
\[ \deg c_\ast(\Xi)=\bb_{3} \quad \iff \quad (P,\Xi)\in \theta^2_{\mathrm{null},4}\,. \] 
where $\theta^2_{\mathrm{null},4}\subset \mathcal{A}_4$ denotes the locus of ppav's with exactly $2$ vanishing theta nulls. By \cite{Debarre1987}, we have
\[ \overline{\theta^2_{\mathrm{null},4}}=\EE_{4,0}\,, \] 
and we treated the case of $\EE_{4,0}$ above already.
\end{proof}

Let $g\geq 4$ and $(P,\Xi)\in \EE_{g,0}^{\prime 0}$ (resp. $(P,\Xi)\in  \SE_g^0\cup \SE_{1,g-1}^0$). Then by Lemma \ref{Lem: Gauss map finite for Eg0'}, Theorem \ref{Theorem: Characteristic Cycle for E'gt} and Lemma \ref{Lem: Chern Mather classes for E'g0} (resp. by \cite[Th. 1,2,3]{podelski2023GaussEgt}):
\begin{itemize}
    \item The Gauss map $\gamma_\Xi: \PP \Lambda_\Xi \to \PP T^\vee_0 P$ is finite away from a locus of codimension at least $3$.
    \item We have
    \[ c_{M,r}(\cc(\IC_\Xi))= \bb_{g-r-1} \frac{[\Xi]^{ g-r}}{(g-r)!} \qquad \text{for $r\geq 1$.}  \]
\end{itemize}
Thus, Theorem \ref{maintheorem: Tannakian Schottky on Bielliptic Prym locus} follows from Theorem \ref{maintheorem: Tannakian criterion for Jacobians}.

\section{Proof of Theorem \ref{maintheorem: Tannakian Schottky in dim up to 5}}\label{Sec: The dimension 5 case}
We now prove that the full Tannakian Schottky conjecture holds in dimension up to $5$. For $g\leq 3$ we have $\mathcal{J}_g=\mathcal{A}_g$ so there is nothing to prove. In dimension $4$, the only missing piece in \cite[Sec. 4]{Kraemer2021MicrolocalGauss2} is the case of non-hyperelliptic fake Jacobians in $\theta^2_{\mathrm{null},4}$, which we treated above already. We now treat the case $g=5$. Recall from (\ref{equ: Andreotti dim 5}) that
\begin{equation*}
    \fJac_5\subset  \mathcal{N}^{(5)}_1 = \mathcal{J}_5 \cup \mathcal{A}_{1,4} \cup \EE_{5,0} \cup \EE_{5,1} \cup \EE_{5,2} \,,
\end{equation*} 
For $(A,\Theta)\in \mathcal{A}_{1,4}$ we have $\chi(\IC_\Theta)=0$. There are no non-hyperelliptic fake Jacobians on $\BE_5=\EE_{5,0} \cup \EE_{5,1} \cup \EE_{5,2}$ by Theorem \ref{maintheorem: Tannakian Schottky on Bielliptic Prym locus}. Suppose from now on that $(A,\Theta)\in \BE_5$ is a hyperelliptic fake Jacobian. We have
\[ G_\Theta= \Sp_{8}(\CC)/ \mu_{2}
       \qquad \text{and}\qquad \omega_\Theta= \Alt^{4}(\CC^{8})/\Alt^{2}(\CC^{8
    })\,. \]
Thus,
 \[ \deg c_\ast(\Theta)=\dim \omega_\Theta =\binom{8}{4}-\binom{8}{2}=42 \,. \] 
 Thus $(A,\Theta)\in \SE_{1,2,2}^2\cup \SE_{2,3}^{10}$ by \cite[Th. 1 and 2]{podelski2023GaussEgt}, \cite[Sec. 2]{podelski2023boundary} and Theorem \ref{Theorem: Characteristic Cycle for E'gt}. By \cite[Th. 4.1.2]{Kraemer2021MicrolocalGauss2} we have 
\[ [2]_\ast \cc(\IC_\Theta)= \mathrm{Alt}^{4}(\Lambda)-\mathrm{Alt}^{2}(\Lambda)\]
for some Lagrangian cycle $\Lambda\in \mathscr{L}(A)$. Thus, by Corollary \ref{Cor: first chern mather classes of Alt^n Lambda} we have in
\begin{equation}\label{Equ: Chern Mather 1 equality contradciting Hyperellitpic Fake Jacobians in dim 5}
    4 c_1(\Theta)=\left(\binom{6}{3}-\binom{5}{1}\right)c_{M,1}(\Lambda)=14 c_{M,1}(\Lambda) \,. 
\end{equation} 
By \cite[Th. 3.7]{podelski2023boundary} we have
\[ \cc(\IC_\Theta)=\Lambda_\Theta+\sum_{\Xiasing} \Lambda_x \,. \]
Let $\ud=(d_1,\dots,d_n)\in \{(1,2,2),(2,2)\}$ and assume $(A,\Theta)\in \SE_\ud$. We use the notation of (\ref{Fig: Cartesian Diagram of curves Egt}). Let $\beta:N''=N_1''\cup\cdots \cup N_n'' \to C''$ be the normalization with $g(N_i)=d_i$. Let $P''=\Prym(C''\to E)$. By \cite{podelski2023boundary} there are $2:1$ isogenies of polarized varieties
\[ P \overset{\pi''^\ast}{\longleftarrow} P'' \overset{\beta^\ast}{\longrightarrow} JN'' \,. \] 
Let $\Xi''=(\pi''^\ast)^\ast \Xi$ and $\theta_1\in \HH^1(P'')$ denote the pullback of the principal polarization of $JN_1$. The conormal varieties to points do not appear in the first Chern-Mather class so
\[ 
c_{M,1}(\cc(\IC_{\Xi''}))\cap \theta_1=c_{M,1}(\Lambda_{\Xi''})\cap \theta_1=
\begin{cases}
    44 &\text{if $\ud=(1,2,2)$,} \\
    92 & \text{if $\ud=(2,3)$,}
\end{cases}
\]
by \cite[Prop. 3.11]{podelski2023boundary}. This is not divisible by $7$ in both cases, contradicting (\ref{Equ: Chern Mather 1 equality contradciting Hyperellitpic Fake Jacobians in dim 5}).

\printbibliography

\end{document}